\documentclass[12pt]{amsart}
\usepackage[]{amscd,amsfonts,amsmath,amsxtra,amssymb}
\usepackage{graphics}
\usepackage[all]{xy}
\newtheorem{sub}{}[section]
\newtheorem{subsub}{}[sub]
\def\ov#1{\overline{#1}}

\def\coker{\mathop{\rm coker}\nolimits}
\def\Hom{\mathop{\rm Hom}\nolimits}
\def\HHom{\mathop{\mathcal Hom}\nolimits}

\def\Ext{\mathop{\rm Ext}\nolimits}
\def\EExt{\mathop{\mathcal Ext}\nolimits}

\def\imm{\mathop{\rm im}\nolimits}
\def\deg{\mathop{\rm deg}\nolimits}
\def\rg{\mathop{\rm rg}\nolimits}
\def\Deg{\mathop{\rm Deg}\nolimits}

\def\lra{\longrightarrow}

\def\sigg{\mathop{\hbox{$\displaystyle\sum$}}\limits}

\def\hfl#1#2{\smash{\mathop{\ \hbox to 12mm{\rightarrowfill}}
\limits^{\scriptstyle#1}_{\scriptstyle#2} \ }}
\def\hflb#1#2{\smash{\mathop{\hbox to 12mm{\leftarrowfill}}
\limits^{\scriptstyle#1}_{\scriptstyle#2}}}

\def\m#1{{\hbox{$#1$}}}

\def\ot{\otimes}
\def\og{\leavevmode\raise.3ex\hbox{$\scriptscriptstyle\langle\!\langle$}}
\def\fg{\leavevmode\raise.3ex\hbox{$\scriptscriptstyle\,\rangle\!\rangle$}}

\def\nsp{\lbrace 0\rbrace}

\def\Ssect#1#2{\pagebreak[3]\begin{sub}\label{#2}{\sc\small\small 
#1}\rm\medskip}

\def\sepsec{\vskip 1.4cm}
\def\sepsub{\vskip 0.7cm}
\def\sepsubsub{\vskip 0.5cm}
\def\sepprop{\vskip 0.5cm}

\def\xmat#1{\[\xymatrix{#1}\]}
\def\flinc{\ar@{^{(}->}}
\def\fleq{\ar@{=}}
\def\flon{\ar@{->>}}
\def\fmaps{\ar@{|-{>}}}
\def\Nligne{\hfil\break}

\newcommand{\M}{{\mathbb M}}

\newcommand{\C}{{\mathbb C}}

\newcommand{\D}{{\mathbb D}}
\newcommand{\F}{{\mathbb F}}
\newcommand{\E}{{\mathbb E}}

\newcommand{\V}{{\mathbb V}}

\renewcommand{\L}{{\mathbb L}}

\def\T{{\mathbb T}}

\newcommand{\ke}{{\mathcal E}}
\newcommand{\kf}{{\mathcal F}}
\newcommand{\kg}{{\mathcal G}}
\newcommand{\kh}{{\mathcal H}}
\newcommand{\ki}{{\mathcal I}}

\newcommand{\km}{{\mathcal M}}
\newcommand{\kn}{{\mathcal N}}
\newcommand{\ko}{{\mathcal O}}

\newcommand{\kt}{{\mathcal T}}
\newcommand{\ku}{{\mathcal U}}

\begin{document}

\def\refname{R\'ef\'erences}
\def\contentsname{Sommaire}
\def\proofname{D\'emonstration}
\def\abstractname{R\'esum\'e}

\author{Jean--Marc Dr\'{e}zet}
\address{
Institut de Math\'{e}matiques de Jussieu\\
UMR 7586 du CNRS\\
175 rue du Chevaleret\newline  
F-75013 Paris, France}
\email{drezet@math.jussieu.fr}

\title[{Conditions d'existence des faisceaux semi-stables}] {Sur les conditions
d'existence des faisceaux semi-stables sur les courbes multiples primitives}

\begin{abstract}
We give sufficient conditions for the (semi-)stability of torsion free sheaves
on a primitive multiple curve. These conditions are used to prove that some
moduli spaces of stable sheaves are not empty. We study mainly the {\em quasi
locally free sheaves of generic type} (this includes the locally free sheaves).
These sheaves are {\em generic}, i.e. for every moduli space of torsion free
sheaves, the sheaves of this type correspond to an open subset of the moduli
space.
\end{abstract}

\maketitle
\tableofcontents

\section{Introduction}
\label{intro}

Une {\em courbe multiple primitive} est une vari\'et\'e alg\'ebrique complexe de
Cohen-Macaulay qui peut localement \^etre plong\'ee dans une surface
lisse, et dont la sous-vari\'et\'e r\'eduite associ\'ee est une courbe lisse.
Les courbes projectives multiples primitives ont \'et\'e d\'efinies et
\'etudi\'ees pour la premi\`ere fois par C.~B\u anic\u a et O.~Forster dans
\cite{ba_fo}. Leur classification a \'et\'e faite dans \cite{ba_ei} pour les
courbes doubles, et dans \cite{dr1} dans le cas g\'en\'eral.
Les faisceaux semi-stables sur des vari\'et\'es non lisses ont d\'ej\`a \'et\'e
\'etudi\'es par diff\'erents auteurs (\cite{ses}, \cite{bho}, \cite{bho2},
\cite{tei1}, \cite{tei2}, \cite{tei3},\cite{in}\cite{in2}).

On peut esp\'erer en trouver des applications concernant les fibr\'es
vectoriels ou leurs vari\'et\'es de modules sur les courbes lisses
(cf. \cite{ei_gr}, \cite{sun0}, \cite{sun1}) en faisant d\'eg\'en\'erer des
courbes lisses vers une courbe multiple primitive. Le probl\`eme de la
d\'eg\'en\'eration des courbes lisses en courbes primitives doubles est
\'evoqu\'e dans \cite{gon}.

Les articles \cite{dr2} et \cite{dr4} sont consacr\'es \`a l'\'etude  des
faisceaux coh\'erents et de leurs vari\'et\'es de modules sur les courbes
multiples primitives. On donne ici des crit\`eres de (semi-)stabilit\'e et des
conditions suffisantes d'existence des faisceaux semi-stables sur ces courbes.
On appliquera ensuite ces crit\`eres \`a des faisceaux sans torsion
g\'en\'eriques. Les conditions existence des faisceaux (semi-)stables
s'expriment en fonction d'{\em invariants} de ces faisceaux, parmi lesquels se
trouvent le rang et le degr\'e g\'en\'eralis\'es.

Le cas des faisceaux localement libres est trait\'e. Dans ce cas les seuls
invariants sont le rang et le degr\'e g\'en\'eralis\'es. Les vari\'et\'es de
modules obtenues sont irr\'eductibles.

On consid\`ere aussi des faisceaux plus compliqu\'es, les {\em faisceaux quasi
localement libres de type rigide} non localement libres, o\`u il y a deux
invariants suppl\'ementaires (dans ce cas les vari\'et\'es de modules de
faisceaux de rang et degr\'e g\'en\'eralis\'es fix\'es peuvent avoir de
multiples composantes).

Pour finir on traitera des exemples simples de faisceaux sans torsion non quasi
localement libres.

\sepsub

\Ssect{Faisceaux coh\'erents sur les courbes multiples primitives}{intro1}

Soient $C$ une courbe projective lisse irr\'eductible, $n$ un entier tel que
\m{n\geq 2} et $Y$ une courbe multiple primitive de multiplicit\'e $n$ et de
courbe r\'eduite associ\'ee $C$. Si \m{\ki_C} est le faisceau d'id\'eaux
de $C$ dans $Y$,
\[L\ = \ \ki_C/\ki_C^2\]
est un fibr\'e en droites sur $C$, dit {\em associ\'e} \`a $Y$. Dans cet article
on supposera que \m{\deg(L)<0}. Le cas o\`u \m{\deg(L)\geq 0} est moins
int\'eressant car les seuls faisceaux stables sont alors les fibr\'es
vectoriels stables sur $C$.

Pour \m{1\leq i\leq n} on note \m{C_i} le sous-sch\'ema de $Y$ d\'efini par le
faisceau d'id\'eaux \m{\ki_C^i}. C'est une courbe multiple primitive de
multiplicit\'e $i$ et on a une filtration
\[C=C_1\subset\cdots\subset C_n=Y\ .\]
On notera \ \m{\ko_i=\ko_{C_i}}.

Le faisceau \m{\ki_C} est localement libre de rang 1 sur \m{C_{n-1}}. Il existe
un fibr\'e en droites $\L$ sur \m{C_n} tel que \m{\L_{\mid C_{n-1}}=\ki_C}. Pour
tout faisceau coh\'erent $\ke$ sur \m{C_n} on a donc un morphisme canonique
\[\ke\ot\L\lra\ke\]
dont le noyau et le conoyau sont ind\'ependants du choix de $\L$.

Si $\kf$ est un faisceau coh\'erent sur $Y$ on
note \m{\kf_i} le noyau de la restriction \ \m{\kf\to\kf_{\mid C_i}},
\m{\kf^{(i)}} celui du morphisme canonique \ \m{\kf\to\kf\ot\L^{-i}}. On a des
suites exactes canoniques
\[0\lra\kf_i\lra\kf\lra\kf_{\mid C_i}\lra 0 ,\]
\[0\lra\kf^{(i)}\lra\kf\lra\kf_i\ot\L^{-i}\lra 0 .\]

Les quotients \ \m{G_i(\kf) \ = \ \kf_i/\kf_{i+1}} , \m{0\leq i<n}, sont des
faisceaux sur $C$. Ils permettent de d\'efinir les {\em rang g\'en\'eralis\'e}
et le {\em degr\'e g\'en\'eralis\'e} de $\kf$ :
\[R(\kf)=\sigg_{i=0}^{n-1}\rg(G_i(\kf)) , \quad\quad
\Deg(\kf)=\sigg_{i=0}^{n-1}\deg(G_i(\kf)) \ .\]
Ce sont des invariants par d\'eformation (cf. \ref{inva}, \cite{dr2},
\cite{dr4}). Si \m{R(\kf)>0}, le nombre rationnel
\[\mu(\kf) \ = \ \frac{\Deg(\kf)}{R(\kf)}\]
s'appelle la {\em pente} de $\kf$.

Pour \m{1\leq i<n}, on note \m{\kf[i]} le noyau du morphisme canonique surjectif
\xmat{\kf\flon[r] & \kf_{\mid C_i}\flon[r] & (\kf_{\mid C_i})^{\vee\vee} .}

\sepprop

\begin{subsub}\label{intro1b}Faisceaux quasi localement libres - \rm
On dit qu'un faisceau coh\'erent $\ke$ sur $Y$ est {\em quasi localement libre}
s'il existe des entiers \m{m_1,\ldots,m_n} non n\'egatifs tels que $\ke$ soit
localement isomorphe \`a
\[\bigoplus_{i=1}^nm_i\ko_{i} .\] 
Les entiers \m{m_i} sont alors uniquement d\'etermin\'es.
\end{subsub}

\sepprop

\begin{subsub}\label{intro1c}Faisceaux quasi localement libres de type rigide -
\rm (cf. \cite{dr4}) Si $\ke$ est quasi localement libre on dit qu'il est {\em
de type rigide} s'il est localement libre, ou s'il existe un unique entier $k$,
\m{1\leq k\leq n-1}, tel que \m{m_k\not=0}, et si on a \m{m_k=1}. Donc un
faisceau quasi localement libre de type rigide non localement libre est
localement isomorphe \`a un faisceau du type \m{a\ko_n\oplus\ko_k}, avec
\m{1\leq k\leq n-1}. L'int\'er\`et de ces faisceaux est que la propri\'et\'e
pour un faisceau d'\^etre quasi localement libre de type rigide est une {\em
propri\'et\'e ouverte}. En particuliers les faisceaux stables localement libre
de type rigide de rang g\'en\'eralis\'e $R$ et de degr\'e g\'en\'eralis\'e $d$
constituent un ouvert de la vari\'et\'e de modules des faisceaux stables de rang
g\'en\'eralis\'e $R$ et de degr\'e g\'en\'eralis\'e $d$ sur \m{C_n}.

Soit $\ke$ un faisceau quasi localement libre de type rigide localement
isomorphe \`a \m{a\ko_n\oplus\ko_k}, avec \m{a\geq 1}, \m{1\leq k<n}. Alors les
faisceaux \m{\ke_k} et \m{\ke^{(k)}} sont localement libres sur \m{C_{n-k}},
\m{C_k} respectivement. On pose
\[E_\ke=\ke_{\mid C} , \quad\quad F_\ke=\ke_{k\mid C} , \ \
V_\ke=(\ke^{(k)})_{\mid C} .\]
Ce sont des fibr\'es vectoriels sur $C$ de rang \m{a+1}, $a$, \m{a+1}
respectivement. On montre en \ref{def_qlltr} qu'on a une suite exacte canonique
\[(*)_\ke \quad\quad\quad 0\lra F_\ke\ot L^{n-k}\lra V_\ke\ot L^{n-k}\lra
E_\ke\lra F_\ke\lra 0 \ .\]
Les rangs et degr\'es g\'en\'eralis\'es des fibr\'es \m{E_\ke} et \m{F_\ke}
(et donc aussi \m{V_\ke}) sont invariants par d\'eformation.
\end{subsub}

\sepprop

\begin{subsub}\label{intro1d}Construction des faisceaux quasi localement libres
de type rigide - \rm Elle est faite par r\'ecurrence sur $n$ dans
\ref{const_000}, \ref{const_qlltr}, \ref{const_qlltr2} et \ref{const_qlltr3}. On
construit le faisceau $\ke$ sur \m{C_n} connaissant \m{\ke_1}, dont le support
est \m{C_{n-1}}, et \m{\ke_{\mid C}}. A priori il semble plus naturel ce
construire $\ke$ connaissant \m{\ke_{\mid C_{n-1}}}. On montre dans
\ref{const_qlltr4} que cela est impossible car les faisceaux sur \m{C_{n-1}} qui
sont des restrictions de faisceaux quasi localement libres de type rigide sur
\m{C_n} sont {\em sp\'eciaux}.

Cette m\'ethode de construction devrait rendre possible la description
pr\'ecise d'ouverts des vari\'et\'es de modules de faisceaux stables qui
contiennent de tels faisceaux.
\end{subsub}

\end{sub}

\sepsub

\Ssect{Vari\'et\'es de modules de faisceaux stables}{intro2}

La (semi-)stabilit\'e au sens de Simpson (cf. \cite{si}) des faisceaux sans
torsion sur \m{C_n} ne d\'epend pas du choix d'un fibr\'e en droites tr\`es
ample sur \m{C_n}. Elle est analogue \`a celle des fibr\'es (semi-)stables sur
les courbes projectives lisses (cf. \cite{dr2}, \cite{dr4}) :
un faisceau sans torsion $\ke$ sur \m{C_n} est {\em semi-stable} (resp.
{\em stable}) si pour tout sous-faisceau propre $\kf$ de $\ke$ on a \
\m{\mu(\kf)\leq\mu(\ke)} (resp. $<$).

L'hypoth\`ese \m{\deg(L)<0} est justifi\'ee par le fait que dans le cas
contraire les seuls faisceaux sans torsion stables sur \m{C_n} sont les fibr\'es
vectoriels stables sur $C$.

Soient $R$, $d$ des entiers, avec \m{R\geq 1}. On note \m{\km(R,d)} la
vari\'et\'e de modules des faisceaux stables de rang g\'en\'eralis\'e $R$ et de
degr\'e g\'en\'eralis\'e $d$ sur \m{C_n}.

Soient $a$, $k$, $\epsilon$, $\delta$ des entiers, avec \m{a\geq 1} et \m{1\leq
k<n}. Soient
\[R \ = \ an+k , \quad\quad\quad d \ = \ k\epsilon+(n-k)\delta+(n(n-1)a+
k(k-1))\frac{\deg(L)}{2} \ .\]
Les faisceaux quasi localement libres $\ke$ de type g\'en\'erique stables
localement isomorphes \`a \m{a\ko_n\oplus\ko_k} et tels que \m{E_\ke} (resp.
\m{F_\ke}) soit de rang \m{a+1} (resp. $a$) et de degr\'e $\epsilon$ (resp.
$\delta$) constituent un ouvert irr\'eductible de \m{\km(R,d)}, dont la
sous-vari\'et\'e r\'eduite sous-jacente est not\'ee
\m{\kn(a,k,\delta,\epsilon)} (cf. \cite{dr4}). A priori \m{\km(R,d)} a donc
plusieurs composantes irr\'eductibles.

\end{sub}

\sepsub

\Ssect{Principaux r\'esultats}{intro3}

On d\'emontre en \ref{theo1} le

\sepprop

{\bf Th\'eor\`eme : }{\em Soient $\ke$ un faisceau
coh\'erent sans torsion sur \m{C_n} et $k$ un entier tel que \m{1\leq k<n} et
que \m{\ke_k\not=0}. On suppose que
\[\mu(\ke^{(k)})\leq\mu(\ke) ,\quad
\mu((\ke^\vee)^{(k)})\leq\mu(\ke^\vee) .\]
Alors, si \m{\ke[k]}, \m{(\ke_{\mid C_k})^{\vee\vee}}, \m{(\ke^\vee)[k]} et
\m{((\ke^\vee)_{\mid C_k})^{\vee\vee}} sont semi-stables il en est de m\^eme de
$\ke$.

Si de plus les in\'egalit\'es de (\ref{eqX}) sont strictes, et si les
faisceaux \m{\ke[k]}, \m{(\ke_{\mid C_k})^{\vee\vee}}, \m{(\ke^\vee)[k]} et
\m{((\ke^\vee)_{\mid C_k})^{\vee\vee}} sont stables, alors $\ke$ est stable.}

\sepprop

M\^eme si on se limitait aux faisceaux quasi localement libres il serait
n\'ecessaire de faire intervenir des sous-faisceaux non quasi localement libres
: on donne en \ref{HN} des exemples de fibr\'es vectoriels sur \m{C_2} dont la
filtration de Harder-Narasimhan comporte des faisceaux non quasi localement
libres.

\sepprop

On applique d'abord le th\'eor\`eme pr\'ec\'edent aux fibr\'es vectoriels :

\sepprop

{\bf Th\'eor\`eme : }{\em Soit $\E$ un fibr\'e vectoriel sur \m{C_n}. Alors,
si \m{\E_{\mid C}} est semi-stable (resp. stable), il en est de m\^eme de $\E$.}

\sepprop

On en d\'eduit que les vari\'et\'es de modules de fibr\'es vectoriels stables
sur \m{C_n} sont non vides, pourvu qu'il n'y ait pas d'incompatibilit\'e au
niveau du rang et du degr\'e g\'en\'eralis\'es. Soient \m{r,\delta} des entiers
avec \m{r\geq 1}. Alors le rang g\'en\'eralis\'e $R$ et le degr\'e
g\'en\'eralis\'e $d$ d'un fibr\'e vectoriel $\E$ sur \m{C_n} tel que
\m{\E_{\mid C}} soit de rang $r$ et de degr\'e $\delta$ sont
\[R=nr , \quad d=n\delta+\frac{n(n-1)}{2}r\deg(L) .\]
L'ouvert \m{U(R,d)} de \m{\km(R,d)} correspondant aux fibr\'es vectoriels
stables est non vide, lisse et irr\'eductible.

\sepprop

On s'int\'eresse ensuite aux faisceaux quasi localement libres de type rigide
non localement libres~:

\sepprop

{\bf Th\'eor\`eme : }{\em Soient $a$, $k$ des entiers tels que \m{a>0} et
\m{1\leq k<n}. Soit $\ke$ un faisceau quasi localement libre de type rigide,
localement isomorphe \`a \m{a\ko_n\oplus\ko_k} et tel que
\[\mu(V_\ke)+\frac{n}{2}\deg(L)\leq\mu(F_\ke)\leq\mu(E_\ke)-\frac{n}{2}\deg(L).
\]
Alors si \m{E_\ke}, \m{F_\ke} et \m{V_\ke} sont semi-stables, il en est de
m\^eme de $\ke$.

Si les in\'egalit\'es pr\'ec\'edentes sont strictes, et si \m{E_\ke}, \m{F_\ke}
et \m{V_\ke} sont stables, il en est de m\`eme de $\ke$.}

\sepprop

Le probl\`eme de l'{\em existence} des faisceaux quasi localement libres de type
rigide (semi-)stables est plus compliqu\'e que celui de l'existence des
fibr\'es vectoriels (semi-)stables, car si $\ke$ en est un, la
(semi-)stabilit\'e de \m{E_\ke}, \m{F_\ke} et \m{V_\ke} impose des conditions
suppl\'ementaires sur les invariants de ces faisceaux, \`a cause de la suite
exacte \m{(*)_\ke}.

Avec les notations de \ref{intro2}, on montrera en \ref{theo4} le

\sepprop

{\bf Th\'eor\`eme : }{\em Si on a
\[\frac{\epsilon}{a+1} \ < \ \frac{\delta}{a} \ < \
\frac{\epsilon-(n-k)\deg(L)}{a+1}\]
\m{\kn(a,k,\delta,\epsilon)} est non vide.}

\sepprop

Ce r\'esultat g\'en\'eralise la proposition 9.2.1 de \cite{dr2}, o\`u le cas
des faisceaux de rang g\'en\'eralis\'e 3 sur \m{C_2} localement isomorphes \`a
\m{\ko_2\oplus\ko_C} \'etait trait\'e. La d\'emonstration du th\'eor\`eme
pr\'ec\'edent utilise la r\'esolution de la {\em conjecture de Lange} (cf.
\cite{rute}).

D'apr\`es \cite{dr4}, prop. 6.4.1, la vari\'et\'e \m{\kn(a,k,\delta,\epsilon)}
est irr\'eductible et lisse, et on a
\[\dim(\kn(a,k,\delta,\epsilon)) \ = \ 1 - \big(\frac{n(n-1)}{2}a^2+k(n-1)a+
\frac{k(k-1)}{2}\big)\deg(L)+(g-1)(na^2+k(2a+1))\]
($g$ d\'esignant le genre de $C$).

\sepprop

On termine par donner des applications du premier des th\'eor\`emes
pr\'ec\'edents \`a des faisceaux non quasi localement libres.

Soient $\E$ un fibr\'e vectoriel sur \m{C_n}, \m{E=\E_{\mid C}} et $Z$ un
ensemble fini de points de $C$. On pose \m{z=h^0(\ko_Z)}. Soient
\m{\phi:\E\to\ko_Z} un morphisme surjectif, et \ \m{\ke_\phi=\ker(\phi)}.
On d\'emontrera en \ref{theo5} le

\sepprop

{\bf Th\'eor\`eme : }{\em Si on a \ \m{z\leq-rg(E)\deg(L)} (resp. $<$) et si $E$
et \m{E_\phi} sont semi-stables (resp. stables), alors $\ke_\phi$ est
semi-stable (resp. stable).}

\end{sub}

\sepsub

\Ssect{Plan des chapitres suivants}{intro4}

Le chapitre 2 contient des rappels sur les courbes multiples primitives et les
propri\'et\'es \'el\'ementaires des faisceaux coh\'erents sur ces courbes. On
d\'ecrit dans \ref{constcoh} la m\'ethode de construction d'un faisceau
coh\'erent $\ke$ sur \m{C_n} connaissant le faisceau \m{\ke_1} sur \m{C_{n-1}}
et \m{\ke_{\mid C}}. Elle sera employ\'ee aussi bien pour les faisceaux
localement libres que pour les faisceaux quasi localement libres de type rigide.
On donne dans \ref{HN} des exemples de fibr\'es vectoriels instables sur une
courbe double primitive dont la filtration de Harder-Narasimhan n'est pas
constitu\'ee de faisceaux quasi localement libres. Cela rend n\'ecessaire, dans
l'\'etude de la (semi-)stabilit\'e d'un faisceau, la consid\'eration de
sous-faisceaux sans torsion g\'en\'eraux dont les filtrations canoniques peuvent
comporter des faisceaux ayant de la torsion.

Le chapitre 3 est une \'etude des faisceaux quasi localement libres de type
rigide et de leur construction.

Le chapitre 4 traite de la dualit\'e des faisceaux coh\'erents sur \m{C_n} et
des faisceaux de torsion.

Dans le chapitre 5 sont d\'emontr\'es les r\'esultats \'enonc\'es dans
\ref{intro3}.
\end{sub}

%\newpage
\sepsec

\section{Pr\'eliminaires}\label{prelim}

\Ssect{D\'efinition des courbes multiples primitives et notations}{def_nota}

Une {\em courbe primitive} est une vari\'et\'e lisse $Y$ de
Cohen-Macaulay telle que la sous-vari\'et\'e r\'eduite associ\'ee \m{C=Y_{red}}
soit une courbe lisse irr\'eductible, et que tout point ferm\'e de
$Y$ poss\`ede un voisinage pouvant \^etre plong\'e dans une surface lisse.

Soient $P$ est un point ferm\'e de $Y$, et $U$ un voisinage de $P$ pouvant
\^etre plong\'e dans une surface lisse $S$. Soit $z$ un \'el\'ement de
l'id\'eal maximal de l'anneau local \m{\ko_{S,P}} de $S$ en $P$ engendrant
l'id\'eal de $C$ dans cet anneau. Il existe alors un unique entier $n$,
ind\'ependant de $P$, tel que l'id\'eal de $Y$ dans \m{\ko_{S,P}} soit
engendr\'e par \m{(z^n)}. Cet entier $n$ s'appelle la {\em multiplicit\'e} de
$Y$. Si \m{n=2} on dit que $Y$ est une {\em courbe double}. D'apr\`es
\cite{dr1}, th\'eor\`eme 5.2.1, l'anneau \m{\ko_{YP}} est isomorphe \`a
\m{\ko_{CP}\ot(\C[t]/(t^n))}.

Soit \m{\ki_C} le faisceau d'id\'eaux de $C$ dans $Y$. Alors le faisceau
conormal de $C$, \m{L=\ki_C/\ki_C^2} est un fibr\'e en droites sur $C$, dit {\em
associ\'e} \`a $Y$. Il existe une filtration canonique
\[C=C_1\subset\cdots\subset C_n=Y\ ,\]
o\`u au voisinage de chaque point $P$ l'id\'eal de \m{C_i} dans \m{\ko_{S,P}}
est \m{(z^i)}. On notera \ \m{\ko_i=\ko_{C_i}} \ pour \m{1\leq i\leq n}.

Le faisceau \m{\ki_C} est un fibr\'e en droites sur \m{C_{n-1}}. Il existe
d'apr\`es \cite{dr2}, th\'eor\`eme 3.1.1, un fibr\'e en droites $\L$ sur
\m{C_n} dont la restriction \`a \m{C_{n-1}} est \m{\ki_C}. On a alors, pour
tout faisceau de \m{\ko_n}-modules $\ke$ un morphisme canonique
\[\ke\ot\L\lra\ke\]
qui en chaque point ferm\'e $P$ de $C$ est la multiplication par $z$.

\end{sub}

\sepsub

\Ssect{Filtrations canoniques}{filtcan}

Dans toute la suite de ce chapitre on consid\`ere une courbe multiple primitive
\m{C_n} de courbe r\'eduite associ\'ee $C$. On utilise les notations de
\ref{def_nota}.

Soient $P$ un point ferm\'e de $C$, $M$ un
\m{\ko_{nP}}-module de type fini. Soit $\ke$ un faisceau coh\'erent sur \m{C_n}.

\sepsubsub

\begin{subsub}\label{QLL-def1}Premi\`ere filtration canonique - \rm
On d\'efinit la {\em premi\`ere filtration canonique de $M$} : c'est la
filtration
\[M_n=\nsp\subset M_{n-1}\subset\cdots\subset M_{1}\subset M_0=M\]
telle que pour \m{0\leq i< n}, \m{M_{i+1}} soit le noyau du morphisme
canonique surjectif \Nligne
\m{M_{i}\to M_{i}\ot_{\ko_{n,P}}\ko_{C,P}} .
On a donc
\[M_{i}/M_{i+1} \ = \ M_{i}\ot_{\ko_{n,P}}\ko_{C,P}, \ \ \ \
M/M_i \ \simeq \ M\ot_{\ko_{n,P}}\ko_{i,P}, \ \ \ \
M_i \ = \ z^iM .\]
 On pose, si $i>0$,
\m{G_i(M)  =  M_i/M_{i+1}} .
Le gradu\'e \ 
\m{{\rm Gr}(M) = \oplus_{i=0}^{n-1}G_i(M) =
\oplus_{i=0}^{n-1}z^iM/z^{i+1}M} \
est un \m{\ko_{C,P}}-module. 

On d\'efinit de m\^eme la {\em premi\`ere filtration canonique de $\ke$} :
c'est la filtration
\[\ke_n=0\subset \ke_{n-1}\subset\cdots\subset \ke_{1}\subset \ke_0=\ke\]
telle que pour \m{0\leq i< n}, \m{\ke_{i+1}} soit le noyau du morphisme
canonique surjectif \ \m{\ke_i\to\ke_{i\mid C}}.
On a donc \
\m{\ke_{i}/\ke_{i+1}=\ke_{i\mid C}} ,
\m{\ke/\ke_i=\ke_{\mid C_i}} .
 On pose, si $i\geq 0$,
\[G_i(\ke) \ = \ \ke_i/\ke_{i+1} .\]
Le gradu\'e
\m{{\rm Gr}(\ke)} est un \m{\ko_{C}}-module.
\end{subsub}

\sepsubsub

\begin{subsub}\label{2-fc}Seconde filtration canonique - \rm
On d\'efinit la {\em seconde filtration canonique de $M$} : c'est la filtration
\[M^{(0)}=\nsp\subset M^{(1)}\subset\cdots\subset M^{(n-1)}\subset M^{(n)}=M\]
avec \
\m{M^{(i)} \ = \ \big\lbrace u\in M ; z^iu=0\big\rbrace} .
Si \ \m{M_n=\nsp\subset M_{n-1}\subset\cdots\subset M_1\subset M_0=M} \ est
la (premi\`ere) filtration canonique de $M$ on a \ \m{M_i\subset M^{(n-i)}} \
pour \m{0\leq i\leq n}. On pose, si $i>0$,
\m{G^{(i)}(M)  =  M^{(i)}/M^{(i-1)}} .
Le gradu\'e \
\m{{\rm Gr}_2(M) = \oplus_{i=1}^nG^{(i)}(M)} \
est un \m{\ko_{C,P}}-module.

On d\'efinit de m\^eme la {\em seconde filtration canonique de $\ke$} :
\[\ke^{(0)}=\nsp\subset \ke^{(1)}\subset\cdots\subset
\ke^{(n-1)}\subset \ke^{(n)}=\ke  .\]
 On pose, si $i>0$,
\[G^{(i)}(\ke) \ = \ \ke^{(i)}/\ke^{(i-1)} .\]
Le gradu\'e
\m{{\rm Gr}_2(\ke)} est un \m{\ko_{C}}-module.
\end{subsub}
\end{sub}

\sepsub

\Ssect{Invariants}{inva}

\begin{subsub}Rang g\'en\'eralis\'e - \rm L'entier \
\m{R(M)=\rg({\rm Gr}(M))} \
s'appelle le {\em rang g\'en\'eralis\'e} de $M$.

L'entier \
\m{R(\ke)=\rg({\rm Gr}(\ke))} \
s'appelle le {\em rang g\'en\'eralis\'e} de $\ke$. On a donc
\m{R(\ke)=R(\ke_P)} pour tout \m{P\in C}.

\end{subsub}

\sepsubsub

\begin{subsub}Degr\'e g\'en\'eralis\'e - \rm L'entier \
\m{\Deg(\ke)=\deg({\rm Gr}(\ke))} \
s'appelle le {\em degr\'e g\'en\'eralis\'e de } $\ke$.

Si \m{R(\ke)>0} on pose \m{\mu(\ke)=\Deg(\ke)/R(\ke)} et on appelle ce nombre la
{\em pente} de $\ke$.
\end{subsub}

\sepprop

Le rang et le degr\'e g\'en\'eralis\'es sont {\em additifs}, c'est-\`a-dire que
si \ \m{0\to\ke'\to\ke\to\ke''\to 0} \ est une suite exacte de faisceaux
coh\'erents sur \m{C_n} alors on a
\[R(\ke)=R(\ke')+R(\ke'') , \quad\quad\quad \Deg(\ke)=\Deg(\ke')+\Deg(\ke'') ,\]
et sont invariants par d\'eformation.

\end{sub}

\sepsub

\Ssect{Faisceaux quasi localement libres}{qll}

Soit $P$ un point ferm\'e de $C$.
Soit $M$ un \m{\ko_{n,P}}-module de type fini.
On dit que $M$ est {\em quasi libre} s'il existe des entiers \m{m_1,\ldots,m_n}
non n\'egatifs et un isomorphisme
\m{M\simeq\bigoplus_{i=1}^n m_i\ko_{i,P}} .
Les entiers \m{m_1,\ldots,m_n} sont uniquement d\'etermin\'es.
On dit que $M$ est {\em de type} \m{(m_1,\ldots,m_n)}. On a \
\m{R(M)=\sigg_{i=1}^ni.m_i}  .

Soit $\ke$ un faisceau coh\'erent sur \m{C_n}.
On dit que $\ke$ est {\em quasi localement libre} en un point $P$ de
$C$ s'il existe un ouvert $U$ de \m{C_n} contenant
$P$ et des entiers non n\'egatifs \m{m_1,\ldots,m_n} tels que pour tout
point $Q$ de $U$, \m{\ke_{n,Q}} soit quasi localement libre de type
\m{m_1,\ldots,m_n}. Les entiers \m{m_1,\ldots,m_n} sont uniquement
d\'etermin\'es et ne d\'ependent que de $\ke$,
et on dit que \m{(m_1,\ldots,m_n)} est le {\em type de } $\ke$. Sur un
voisinage de $P$, $\ke$ est alors isomorphe \`a \m{\oplus_{i=1}^nm_i\ko_i}.

On dit que $\ke$ est {\em quasi localement libre} s'il l'est en tout point de
\m{C_n}.

On montre que $\ke$ est quasi localement libre en $P$ si et seulement si
pour \m{0\leq i<n}, \m{G_i(\ke)} est libre en $P$.

Il en d\'ecoule que $\ke$ est quasi localement libre si et seulement si
pour \m{0\leq i<n}, \m{G_i(\ke)} est localement libre sur $C$.

\end{sub}

\sepsub

\Ssect{Construction des faisceaux coh\'erents}{constcoh}

\begin{subsub}\label{constcoh1}\rm On d\'ecrit ici le
moyen de construire un faisceau coh\'erent $\ke$ sur \m{C_n}, connaissant
\m{\ke_{\mid C}} et \m{\ke_1}, qui sont des faisceaux sur $C$ et \m{C_{n-1}}
respectivement.

Soient $\kf$ un faisceau coh\'erent sur \m{C_{n-1}} et $E$ un fibr\'e vectoriel
sur $C$. On s'int\'eresse aux faisceaux coh\'erents $\ke$ sur \m{C_n} tels que
\m{\ke_{\mid C}=E} et \m{\ke_1=\kf}. Soit \
\m{0\to\kf\to\ke\to E\to 0} \ une suite exacte, associ\'ee \`a \
\m{\sigma\in\Ext_{\ko_n}^1(E,\kf)}.
On voit ais\'ement que le morphisme canonique \ \m{\pi_\ke:\ke\ot\ki_C\to\ke} \
induit un morphisme
\[\Phi_{\kf,E}(\sigma) :E\ot L\lra\kf_{\mid C} .\]
On a \m{\ke_{\mid C}=E} et \m{\ke_1=\kf} si et seulement si
\m{\Phi_{\kf,E}(\sigma)} est surjectif (\cite{dr4}, lemme 3.5.2).

On a une suite exacte canonique
\begin{equation}\label{equCC}
\xymatrix{0\ar[r] & \Ext^1_{\ko_C}(E,\kf^{(1)})\ar[r] & \Ext^1_{\ko_n}(E,\kf)
\ar[rr]^-{\Phi_{\kf,E}} & & \Hom(E\ot L,\kf_{\mid C})\ar[r] & 0 .}
\end{equation}
(\cite{dr4}, prop. 3.5.3).
\end{subsub}

\sepprop

\begin{subsub}\label{constcoh2}\rm On suppose que \m{n\geq 3}. On s'int\'eresse
maintenant aux extensions \ \ \m{0\to\kf\to\ke\to E\to 0} \ associ\'ees aux \
\m{\sigma\in\Ext_{\ko_n}^1(E,\kf)} \ tels que \ \m{\Phi_{\kf,E}(\sigma)=0}
(donc \m{\sigma\in\Ext^1_{\ko_C}(E,\kf^{(1)})} ). Dans ce cas $\ke$ est
localement isomorphe \`a \m{\kf\oplus E} (cf. \cite{dr4}, 2.4). Plus
pr\'ecis\'ement dans la suite exacte (\ref{equCC}) le terme
\m{\Ext^1_{\ko_C}(E,\kf^{(1)})} est en fait \m{H^1(\HHom(E,\kf))}. On peut donc
repr\'esenter $\sigma$ par un cocycle \m{(f_{ij})} relativement \`a un
recouvrement ouvert \m{(U_i)} de \m{C_n}, \m{f_{ij}} \'etant un morphisme \
\m{E_{\mid U_{ij}}\to\kf_{\mid U_{ij}}}. Le faisceau $\ke$ est obtenu en
recollant les \m{(\kf\oplus E)_{\mid U_i}} au moyen des morphismes
\m{\begin{pmatrix} I_\kf & f_{ij}\\ 0 & I_E\end{pmatrix}} (\cite{dr4}, prop.
2.4.3).

On suppose maintenant que $\kf$ est localement libre sur \m{C_{n-1}}. Soit
\m{F=\kf_{\mid C}\ot L^{-1}}, on a donc \ \m{\kf^{(1)}=F\ot L^{n-1}}. En
utilisant la construction pr\'ec\'ente de $\ke$ au moyen d'un cocycle on voit
ais\'ement que \ \m{\ke_{\mid C}\simeq(F\ot L)\oplus E}, et qu'on a une suite
exacte
\[0\lra F\ot L^{n-1}\lra\ke^{(1)}\lra E\lra 0 ,\]
qui est associ\'ee \`a $\sigma$.
\end{subsub}

\sepprop

\begin{subsub}\label{constcoh3} Construction des fibr\'es vectoriels - \rm
On suppose que $\kf$ est un fibr\'e vectoriel sur \m{C_{n-1}}. On veut
construire et param\'etrer les fibr\'es vectoriels $\E$ sur \m{C_n} tels que
\m{\E_1=\kf}. Il convient donc de prendre \ \m{E=\kf_{\mid C}\ot L^{-1}} et de
consid\'erer les extensions \ \m{0\to\kf\to\E\to E\to 0} \ telles que
l'\'el\'ement associ\'e $\sigma$ de \m{\Ext_{\ko_n}^1(E,\kf)} soit tel que \
\m{\Phi_{\kf,E}(\sigma):E\ot L\to E\ot L} \ soit l'identit\'e de \m{E\ot L}.
Si $E$ est simple on montre ais\'ement, en utilisant le fait que \m{\deg(L)<0
}, que deux \'el\'ements \m{\sigma,\sigma'} de \m{\Phi_{\kf,E}^{-1}(I_{E\ot L})}
d\'efinissent des fibr\'es vectoriels $\E$ isomorphes si et seulement si
\m{\sigma=\sigma'}. Dans ce cas les fibr\'es vectoriels recherch\'es sont donc
param\'etr\'es par l'espace affine \
\m{\Phi_{\kf,E}^{-1}(I_E)\simeq\Ext^1_{\ko_C}(E,E\ot L^{n-1})}.
\end{subsub}

\end{sub}

\sepsub

\Ssect{Filtration de Harder-Narasimhan}{HN}

Rappelons qu'on suppose que \m{\deg(L)<0}.
On montre ici que la filtration de Harder-Narasimhan d'un fibr\'e vectoriel sur
\m{C_2} n'est pas n\'ecessairement constitu\'ee de faisceaux quasi localement
libres. Cela entraine que dans l'\'etude de la (semi-)stabilit\'e des faisceaux
localement libres (ou a fortiori quasi localement libres) il faut aussi
consid\'erer des sous-faisceaux sans torsion non n\'ecessairement quasi
localement libres.

Soient $P$ un point ferm\'e de \m{C_2} et \m{\ki_P} son faisceau d'id\'eaux.
Soient \m{z\in\ko_{2,P}} un g\'en\'erateur de l'id\'eal de $C$ et
\m{x\in\ko_{2,P}} au dessus d'un g\'en\'erateur de l'id\'eal de $P$ dans
\m{\ko_{C,P}}. On a donc \m{\ki_{P,P}=(x,z)}. On a une suite exacte de
\m{\ko_{2,P}}-modules
\begin{equation}\label{eqm}
\xymatrix{0\ar[r] & (x,z)\ar[r]^-\alpha & 2\ko_{2,P}\ar[r]^-\beta & (x,z)\ar[r]
& 0 \ ,}\end{equation}
o\`u pour tous \m{a,b\in\ko_{2,P}}
\[\alpha(ax+bz)=(-az,ax+bz) , \quad\quad \beta(a,b)=ax+bz \ .\]
On va globaliser cette suite exacte afin d'obtenir des suites exactes
\begin{equation}\label{eqm0}0\lra\ki_P\ot\D\lra\E\lra\ki_P\lra 0 \ ,
\end{equation}
o\`u $\D$ est fibr\'e en droites sur \m{C_2} et $\E$ un fibr\'e vectoriel de
rang 2 sur \m{C_2}. Le faisceau \m{\EExt^1_{\ko_2}(\ki_P,\ki_P\ot\D)} est
concentr\'e au point $P$.
%On montre ais\'ement qu'au point $P$ ce faisceau est
%isomorphe au \m{\ko_{2,P}}-module \m{2\C}, en utilisant la r\'esolution libre
%\xmat{\cdots\ar[r] & 2\ko_{2,P}\ar[r]^f & 2\ko_{2,P}\ar[r]^g &
%2\ko_{2,P}\ar[r]^q & (x,z) ,}
%o\`u \m{f=\begin{pmatrix}z & 0\\ x & z\end{pmatrix}}, \m{g=\begin{pmatrix}z &
%0\\ -x & z\end{pmatrix}}, \m{q=(x,z)}. 
On en d\'eduit qu'il existe une section
$s$ de \m{\EExt^1_{\ko_2}(\ki_P,\ki_P\ot\D)} dont la valeur en $P$ correspond
\`a l'extension (\ref{eqm}).

On a un morphisme surjectif canonique
\[\Psi:\Ext^1_{\ko_2}(\ki_P,\ki_P\ot\D)\lra
H^0(\EExt^1_{\ko_2}(\ki_P,\ki_P\ot\D)) \ .\]
Donc \m{\Psi^{-1}(s)} est non vide. Si \ \m{0\to\ki_P\ot\D\to\ke\to\ki_P\to 0} \
est une extension associ\'ee à un \'el\'ement de \m{\Psi^{-1}(s)}, le faisceau
$\ke$ est localement libre. L'existence des extensions (\ref{eqm0}) est donc
prouv\'ee.

\sepprop

\begin{subsub}\label{proHN}{\bf Proposition : } Soit \
\m{0\to\ki_P\ot\D\to\E\to\ki_P\to 0} \ une extension, o\`u $\D$ est un fibr\'e
en droites sur \m{C_2} et $\E$ un fibr\'e vectoriel de rang 2 sur \m{C_2}. Alors
si \ \m{\deg(\E_{\mid C})>0} , le faisceau \m{\ki_P\ot\D} est le sous-faisceau
semi-stable maximal de $\E$.
\end{subsub}

\begin{proof} Soit \m{\kh\subset\E} le sous-faisceau semi-stable maximal de
$\E$. On a \ \m{R(\kh)=1,2} ou $3$.

On note $\L_x$ le faisceau d'id\'eaux \'egal \`a \m{\ko_2} sur
\m{C_2\backslash P} et \`a \m{(x)} au point $P$. C'est un fibr\'e en droites
sur \m{C_2} et on a \
\m{\ki_P^\vee\simeq\ki_P\ot\L_x^{-1}} .
On a donc une suite exacte
\[0\lra\ki_P\ot\D\lra\E^\vee\ot\L_x\ot\D\lra\ki_P\lra 0 \]
En consid\'erant cette suite exacte on se ram\`ene au cas o\`u \m{R(\kh)=1} ou
$2$.
 
Soit \m{\kf\subset\ki_P} un sous-faisceau propre. On a alors
\m{R(\kf)=1}, donc $\kf$ est concentr\'e sur $C$, et est donc contenu dans
\m{(\ki_P)^{(1)}=L}. Donc
\[\mu(\kf) \ \leq \ \deg(L) \ \leq \ \mu(\ki_P)=\frac{\deg(L)-1}{2} \ ,\]
(car \m{\deg(L)<0}).

Supposons d'abord que \m{R(\kh)=1}. Si \ \m{\kh\subset\ki_P\ot\D} \ on a
\ \m{\mu(\kh)\leq\mu(\ki_P\ot\D)} ,
ce qui contredit la maximalit\'e de $\kh$. Si \ \m{\kh\not\subset\ki_P\ot\D}, on
peut voir $\kh$ comme un sous-faisceau de \m{\ki_P}, donc \
\m{\mu(\kh)\leq\mu(\ki_P)} ,\ 
donc \m{\mu(\kh)>\mu(\ki_P\ot\D)}, ce qui est absurde.

On a donc \m{R(\kh)=2}. Soit $r$ le rang g\'en\'eralis\'e de l'image $\ku$ de
$\kh$ dans \m{\ki_P}. Si \m{r=0} on a \m{\kh=\ki_P\ot\D}, ce qu'il fallait
d\'emontrer. Si \m{r=2} on peut voir $\kh$ comme un sous-faisceau de \m{\ki_P}
et on a alors encore \
\m{\mu(\kh)\leq\mu(\ki_P)} ,
ce qui est impossible.

Il reste \`a traiter le cas o\`u \m{r=1} en montrant qu'il est impossible. Soit
$d$ le degr\'e de $\ku$ (qui est concentr\'e sur $C$). On a, puisque \
\m{\kh\cap(\ki_P\ot\D)} \ est aussi de rang g\'en\'eralis\'e 1
\[\deg(\ku)\leq\deg(L)\quad\text{et}\quad
\deg(\kh\cap(\ki_P\ot\D))\leq\deg(L)+\deg(\D_{\mid C}) .\]
Donc
\[\mu(\kh) \ \leq \ \deg(L)+\frac{\deg(\D_{\mid C})}{2} \ < \
\frac{\deg(L)-1}{2}+\deg(\D_{\mid C})=\mu(\ki_P\ot\D) ,\]
ce qui contredit la d\'efinition de $\kh$.
\end{proof}

\sepprop

\begin{subsub}\label{rem}{\bf Remarque : }\rm Si on suppose que \m{\deg(\D_{\mid
C})=0} on obtient des fibr\'es vectoriels $\E$ semi-stables de rang 2 sur
\m{C_2} dont la filtration de Jordan-H\"older n'est pas constitu\'ee de
faisceaux quasi localement libres.
\end{subsub}

\end{sub}

\sepsec
%\newpage

\section{Faisceaux quasi localement libres de type rigide}\label{qlltr}

Dans toute la suite de ce chapitre on consid\`ere une courbe multiple primitive
\m{C_n} de courbe r\'eduite associ\'ee $C$. On utilise les notations de
\ref{def_nota}, et on suppose que \m{\deg(L)<0}.

\sepsub

\Ssect{D\'efinitions}{def_qlltr}

Soit $\ke$ un faisceau coh\'erent quasi localement libre sur \m{C_n}. Soient
\m{a=\lbrack\frac{R(\ke)}{n}\rbrack} et \m{k=R(\ke)-an}. On a donc
\m{R(\ke)=an+k}. On dit que $\ke$ est de {\em type rigide} s'il est localement
libre si \m{k=0}, et localement isomorphe \`a \m{a\ko_n\oplus\ko_k} si \m{k>0}.
Si \m{k>0} cela revient à dire que $\ke$ est de type \m{(m_1,\ldots,m_n)}, avec
\m{m_i=0} si \m{i\not=k,n} et \m{m_k=0} ou $1$.

Le fait d'\^etre quasi localement libre de type rigide est une {\em
propri\'et\'e ouverte} : autrement dit si $Y$ une vari\'et\'e
alg\'ebrique int\`egre et $\kf$ une famille plate de faisceaux coh\'erents sur
\m{C_n} param\'etr\'ee par $Y$, alors l'ensemble des points \m{y\in Y} tels que
\m{\ke_y} soit quasi localement libre de type rigide est un ouvert de $Y$
(\cite{dr4}, prop. 6.3.1).

Supposons que $\ke$ soit quasi localement libre de type rigide et que \m{k>0}.
Alors \ \m{\E=\ke_{\mid C_k}} \ est un fibr\'e vectoriel de rang \m{a+1} sur
\m{C_k}, et \ \m{\F=\ke_k} \ est un fibr\'e vectoriel de rang $a$ sur
\m{C_{n-k}}. Donc $\ke$ est une extension
\[0\lra\F\lra\ke\lra\E\lra 0\]
d'un fibr\'e vectoriel de rang \m{a+1} sur \m{C_k} par un fibr\'e vectoriel de
rang $a$ sur \m{C_{n-k}}. De m\^eme \ \m{\V=\ke^{(k)}} \ est un fibr\'e
vectoriel sur \m{C_k} et on a une suite exacte
\[0\lra\V\lra\ke\lra\F\ot\L^{-k}\lra 0 \ .\]

Posons \
\m{E=\ke_{\mid C}=\E_{\mid C}} , \m{F=G_k(\ke)\ot L^{-k}=\F_{\mid C}\ot
L^{-k}} .
Alors on a \ \m{\rg(E)=a+1} , \m{\rg(F)=a}, et
\[\big(G_0(\ke),G_1(\ke),\ldots,G_{n-1}(\ke)\big) \ = \ (E,E\ot L,\ldots,E\ot
L^{k-1},F\ot L^k,\ldots, F\ot L^{n-1}) \ .\]
Donc
\[\Deg(\ke) \ = \ k\deg(E)+(n-k)\deg(F)+\big(n(n-1)a+
k(k-1)\big)\deg(L)/2 \ .\]
On a
\[G^{(n)}(\ke)=\ke/\ke^{(n-1)}=\ke_{n-1}\ot L^{1-n}=G_{n-1}(\ke)\ot L^{1-n}=F
\ .\]
Posons
\m{V=G^{(k)}(\ke)\ot L^{k-n}=\V_{\mid C}\ot L^{k-n}} .
On a
\m{\rg(V)=a+1}, \Nligne\m{\deg(V)=\deg(E)-(n-k)\deg(L)} ,
et
\[\big(G^{(n)}(\ke),G^{(n-1)}(\ke),\ldots,G^{(1)}(\ke)\big) \ = \
(F,F\ot L,\ldots,F\ot L^{n-k-1},V\ot L^{n-k},\ldots,V\ot L^{n-1}\big) \ .\]
Les morphismes canoniques
\xmat{G_i(\ke)\ot L\flinc[r] & G_{i+1}(\ke) \ , & G^{(i+1)}\ot L\flon[r] &
G^{(i)}(\ke)}
d\'efinissent un morphisme surjectif \ \m{\phi:E\to F} \ et un morphisme
injectif \ \m{\psi:F\to V} . D'apr\`es \cite{dr4}, cor. 3.1.8, on a un
isomorphisme canonique
\[\ker(\phi) \ \simeq \ \coker(\psi)\ot L^{n-k} .\]
Posons \ \m{D=\ker(\phi)} . C'est un fibr\'e en droites sur $C$. On a des
suites exactes
\xmat{0\ar[r] & D\ar[r] & E\ar[r]^\phi & F\ar[r] & 0 ,}
\xmat{0\ar[r] & F\ar[r]^\psi & V\ar[r] & D\ot L^{k-n}\ar[r] & 0 .}

\sepprop

\begin{subsub}{\bf Notations : }\rm On pose \ \m{E_\ke=E}, \m{F_\ke=F},
\m{V_\ke=V},
\[\phi_\ke=\phi : E_\ke\lra F_\ke ,\quad\quad
\psi_\ke=\phi : F_\ke\lra V_\ke ,\]
et \ \m{D_\ke=D}. On a une suite exacte canonique
\[(*)_\ke \quad\quad\quad 0\lra F_\ke\ot L^{n-k}\lra V_\ke\ot L^{n-k}\lra
E_\ke\lra F_\ke\lra 0 \ .\]
\end{subsub}

\sepprop

\begin{subsub}\label{const_000} Construction et param\'etrisation - \rm
On cherche ici \`a d\'ecrire comment on peut obtenir les faisceaux quasi
localement libres de type rigide $\ke$ pr\'ec\'edents. On part d'abord d'un
fibr\'e vectoriel $\F$ sur \m{C_{n-k}} de rang \m{a\geq 1} (cf. \ref{constcoh3}
pour la construction et la param\'etrisation des fibr\'es vectoriels) qui sera
\m{\ke_k}. On construira ensuite successivement \m{\ke_{k-1},\cdots,\ke_1,\ke}.
Il y a deux cas diff\'erents : le passage de $\F$ \`a \m{\ke_{k-1}}, et celui
de \m{\ke_i} \`a \m{\ke_{i-1}} si \m{1\leq i<k}. On va donc \'etudier dans les
sections suivantes les deux \'etapes suivantes :

La premi\`ere \'etape consiste \`a \'etudier les extensions
\[0\lra\F\lra\ke\lra H\lra 0\]
sur \m{C_{n-k+1}}, o\`u $H$ est un fibr\'e vectoriel de rang \m{a+1} sur $C$,
telles que le morphisme induit \m{\Phi_{\F,H}:H\to\F_{\mid C}} soit surjectif
(cf. \ref{constcoh}). Le faisceau $\ke$ est alors quasi localement libre de type
rigide, et localement isomorphe \`a \m{a\ko_{n-k+1}\oplus\ko_C} . On a
\m{\ke_{\mid C}=H} et \m{\ke_1=\F}.

Dans la seconde \'etape on part d'un faisceau quasi localement libre de type
rigide $\kg$ sur \m{C_{n-k+i}}, \m{1\leq i<k}, localement isomorphe \`a
\m{a\ko_{n-k+i}\oplus\ko_i}. Soit \m{H=\kg_{\C}\ot L^{-1}}. On s'int\'eresse
alors aux extensions
\[0\lra\kg\lra\ke\lra H\lra 0\]
sur \m{C_{n-k+i+1}} telles que le morphisme induit \m{\Phi_{\kg,H}:H\ot L\to
H\ot L} soit l'identit\'e de \m{H\ot L}. Le faisceau $\ke$ est alors quasi
localement libre de type rigide, et localement isomorphe \`a
\m{a\ko_{n-k+i+1}\oplus\ko_{i+1}}. On a \m{\ke_{\mid C}=H} et \m{\ke_1=\kg}.
\end{subsub}

\end{sub}

\sepsub

\Ssect{Construction et param\'etrisation - Premi\`ere \'etape}{const_qlltr}

On d\'ecrit ici la premi\`ere \'etape \'evoqu\'ee dans \ref{const_000}, dont on
conserve les notations.

\sepprop

On pose \m{F=\F_{\mid C}\ot L^{-1}}. Soient \
\m{\sigma\in\Ext^1_{\ko_{n-k}}(H,\F)} \ et
\[0\lra\F\lra\ke_\sigma\lra H\lra 0\]
l'extension correspondante. On suppose que \
\m{\phi=\Phi_{\F,H}(\sigma)\ot I_{L^-1}:H\to F} \ est surjectif. Soit
\m{D=\ker(\phi)}. On a \m{E_{\ke_\sigma}=H}, \m{F_{\ke_\sigma}=F}, et une suite
exacte
\begin{equation}\label{equDD}
0\lra F\ot L^{n-k}\lra V_{\ke_\sigma}\ot L^{n-k}\lra D\lra 0 .
\end{equation}
On a d'apr\`es \ref{constcoh1} une suite exacte
\xmat{0\ar[r] & \Ext^1_{\ko_C}(H,F\ot L^{n-k})\ar[r] &
\Ext^1_{\ko_{n-k+1}}(H,\F)\ar[rr]^-{\Phi_{\F,H}} & & \Hom(H,F)\ar[r] & 0 .}

\sepprop

\begin{subsub}\label{lemC}
{\bf Lemme : } L'image de $\sigma$ dans \m{\Ext^1_{\ko_{n-k+1}}(D,\F)} est
contenue dans \m{\Ext^1_{\ko_C}(D,F\ot L^{n-k})}, et c'est l'\'el\'ement
associ\'e \`a la suite exacte (\ref{equDD}).
\end{subsub}

\begin{proof}
Soit \m{\sigma'} l'image de $\sigma$ dans \m{\Ext^1_{\ko_{n-k+1}}(D,\F)}.
Notons que la fonctorialit\'e de \m{\phi_{\F,H}} par rapport \`a $\F$ et $H$
entraine que \m{\Phi_{\F,D}(\sigma')=0}. On a donc bien d'apr\`es
\ref{constcoh2} \m{\sigma'\in\Ext^1_{\ko_C}(D,F\ot L^{n-k})}. D'apr\`es
\cite{dr3}, prop. 4.3.1 on a un diagramme commutatif avec lignes exactes
\xmat{
0\ar[r] & \F\ar[r]\fleq[d] & \kh\ar[r]\flinc[d] & D\ar[r]\flinc[d] & 0\\
0\ar[r] & \F\ar[r] & \ke_\sigma\ar[r] & H\ar[r] & 0}
l'extension du haut \'etant associ\'ee \`a \m{\sigma'}. On a \
\m{\kh^{(1)}\subset\ke_\sigma^{(1)}}, et d'apr\`es \ref{constcoh2} on a une
suite exacte
\[0\lra F\ot L^{n-k}\lra \kh^{(1)}\ot L^{n-k}\lra D\lra 0 .\]
Il en d\'ecoule que \ \m{\kh^{(1)}=\ke_\sigma^{(1)}=V_{\ke_\sigma}}. D'apr\`es
\ref{constcoh2} \m{\sigma'} correspond bien \`a l'extension (\ref{equDD}).
\end{proof}

\sepprop

\begin{subsub}\label{propC}{\bf Proposition : } Pour toute extension
\[0\lra F\ot L^{n-k}\lra W\ot L^{n-k}\lra D\lra 0\]
sur $C$ il existe \ \m{\sigma_0\in\Phi_{H,\F}^{-1}(\phi\ot I_L)} \ tel que
l'extension
pr\'ec\'edente soit isomorphe \`a l'extension
\[0\lra F\ot L^{n-k}\lra V_{\ke_{\sigma_0}}\ot L^{n-k}\lra D\lra 0\]
\end{subsub}

\begin{proof}
Cela d\'ecoule du lemme \ref{lemC}, du carr\'e commutatif
\xmat{
\Ext^1_{\ko_C}(H,F\ot L^{n-k})\flinc[r]\ar[d] &
\Ext^1_{\ko_{n-k+1}}(H,\F)\ar[d] \\
\Ext^1_{\ko_C}(D,F\ot L^{n-k})\flinc[r] & \Ext^1_{\ko_{n-k+1}}(D,\F)
}
et de la surjectivit\'e du morphisme de gauche.
\end{proof}

\sepprop

Soient \m{\phi:H\to\F_{\mid C}} un morphisme surjectif et \
\m{\eta\in\Ext^1_{\ko_{n-k+1}}(D,\F)}. Alors on a\Nligne
\m{\Phi_{\F,H}(\sigma)=\phi\ot I_L} et la suite exacte
\[0\lra F\ot L^{n-k}\lra V_{\ke_{\sigma_0}}\ot L^{n-k}\lra D\lra 0\]
est associ\'ee \`a $\eta$ si et seulement si $\eta$ appartient au sous-espace
affine \ \m{\Phi_{\F,H}^{-1}(\phi\ot I_L)\cap\psi^{-1}(\eta)} \ de
\m{\Ext^1_{\ko_{n-k+1}}(H,\F)} ($\psi$ d\'esignant l'application canonique \
\m{\Ext^1_{\ko_{n-k+1}}(H,\F)\to\Ext^1_{\ko_{n-k+1}}(D,\F)} ).

\end{sub}

\sepsub

\Ssect{Construction et param\'etrisation - Seconde \'etape}{const_qlltr2}

On d\'ecrit ici la seconde \'etape \'evoqu\'ee dans \ref{const_000}, dont on
conserve les notations.

On suppose que \m{H=\kg_{\mid C}\ot L^{-1}}. Soient \
\m{\sigma\in\Ext^1_{\ko_{n-k+i+1}}(H,\kg)} \ tel que
\m{\Phi_{\kg,H}(\sigma)} soit l'identit\'e de \m{H\ot L} et \
\m{0\to\kg\to\ke_\sigma\to H\to 0} \ l'extension correspondante.

\sepprop

\begin{subsub}\label{propC2}{\bf Proposition : } On a \m{E_{\ke_\sigma}=E_\kg\ot
L^{-1}}, \m{F_{\ke_\sigma}=F_\kg\ot L^{-1}}, \m{V_{\ke_\sigma}=V_\kg\ot L^{-1}}
et\Nligne \m{(*)_{\ke_\sigma}=(*)_\kg\ot L^{-1}}.
\end{subsub}

\begin{proof} Il suffit de le faire avec \m{a\ko_{n-k+i+1}\oplus\ko_{i+1}} \`a
la place de \m{\ke_\sigma} en utilisant les isomorphismes
locaux \m{\ke_\sigma\simeq a\ko_{n-k+i+1}\oplus\ko_{i+1}} et la fonctorialit\'e
de \m{(*)_{\ke_\sigma}}, ce qui est imm\'ediat.
\end{proof}

\sepprop

On a d'apr\`es \ref{constcoh1} une suite exacte
\xmat{0\ar[r] & \Ext^1_{\ko_C}(H,\kg^{(1)})\ar[r] &
\Ext^1_{\ko_{n-k+i+1}}(H,\kg)
\ar[rr]^-{\Phi_{\kg,H}} & & \Hom(H\ot L,\kg_{\mid C})\ar[r] & 0 \ .}
Les faisceaux \m{\ke_\sigma} consid\'er\'es ici sont donc index\'es par le
sous-espace affine \m{\Phi_{\kg,H}^{-1}(I_{H\ot L})} de
\m{\Ext^1_{\ko_{n-k+i+1}}(H,\kg)}.

\end{sub}

\sepsub

\Ssect{Construction et param\'etrisation - Conclusion}{const_qlltr3}

\begin{subsub}\label{propC3}{\bf Proposition : } Soient $k$, $a$ des
entiers tels que \m{1\leq k<n}, \m{a>0}. Soient $E$, $F$, $V$ des fibr\'es
vectoriels sur $C$ de rangs \m{a+1}, $a$, \m{a+1} respectivement, et
\begin{equation}\label{eqCC2}
0\lra F\ot L^{n-k}\lra V\ot L^{n-k}\lra E\lra F\lra 0
\end{equation}
une suite exacte. Alors il existe un faisceau quasi localement libre de type
rigide $\ke$, localement isomorphe \`a \m{a\ko_n\oplus\ko_k} et tel que
\m{(*)_\ke} soit isomorphe \`a (\ref{eqCC2}).
\end{subsub}

\sepprop

Cela signifie qu'il existe un diagramme commutatif reliant les suite exactes
\m{(*)_\ke} et (\ref{eqCC2}) :
\xmat{F\ot L^{n-k}\ar[r]\ar[d]^\simeq & V\ot L^{n-k}\ar[r]\ar[d]^\simeq &
E\ar[r]\ar[d]^\simeq & F\ar[d]^\simeq \\
F_\ke\ot L^{n-k}\ar[r] & V_\ke\ot L^{n-k}\ar[r] & E_\ke\ar[r] &  F_\ke }
\end{sub}

\sepsub

\Ssect{Restrictions des faisceaux quasi localement libres de type rigide}
{const_qlltr4}

Les m\'ethodes pr\'ec\'edentes de construction de faisceaux quasi localement
libres de type rigide se font sur le principe suivant : on part d'un tel
faisceau $\kf$ sur \m{C_{n-1}} et on en construit un $\ke$ sur \m{C_n} tel que
\m{\ke_1=\kf}.

A priori il semblerait plus naturel de chercher un faisceau $\ke$ tel que \
\m{\ke_{\mid C_{n-1}}=\kf}. Mais c'est impossible car un faisceau quasi
localement libre de type rigide sur \m{C_{n-1}}, non localement libre, n'est pas
n\'ecessairement la restriction d'un faisceau du m\`eme type sur \m{C_n} :

\sepprop

\begin{subsub}\label{propSP}{\bf Proposition : } Soit $\ke$ un faisceau quasi
localement libre de type rigide non localement libre sur \m{C_n} localement
isomorphe \`a \m{\ko_n\oplus\ko_k}, avec \m{a\geq 1} et \m{1\leq k<n-1}. Alors
\m{(\ke_{\mid C_{n-1}})^{(1)}} est scind\'e.
\end{subsub}

\begin{proof} Soient $P$ un point ferm\'e de \m{C_n} et \m{z\in\ko_{n,P}} un
g\'en\'erateur de l'id\'eal de $C$. On fixe un isomorphisme \ \m{\ke_P\simeq
a\ko_{n,P}\oplus\ko_{k,P}}. On a alors \ \m{(\ke_{\mid
C_{n-1}})_P=a\ko_{n-1,P}\oplus\ko_{k,P}} , et
\[(\ke^{(1)})_P=a(z^{n-1})\oplus(z^{k-1}) , \quad
\big((\ke_{\mid C_{n-1}})^{(1)}\big)_P=a(z^{n-2})/(z^{n-1})\oplus(z^{k-1}) .\]
L'image du morphisme canonique \ \m{\lambda:\ke^{(1)}\to(\ke_{\mid
C_{n-1}})^{(1)}} \ au point $P$ est \m{(z^{k-1})}. L'autre facteur
\m{a(z^{n-2})/(z^{n-1})} est
\m{\big((\ke_{\mid C_{n-1}})_{n-2}\big)_P}. On a donc
\[(\ke_{\mid C_{n-1}})^{(1)} \ = \ \imm(\lambda)\oplus(\ke_{\mid C_{n-1}})_{n-2}
.\]
\end{proof}

\end{sub}

\sepsec
%\newpage

\section{Dualit\'e et torsion}\label{f_to}

On consid\`ere dans ce chapitre une courbe multiple primitive \m{C_n} de courbe
r\'eduite associ\'ee $C$. On utilise les notations de \ref{def_nota}.

\sepsub

\Ssect{G\'en\'eralit\'es sur la dualit\'e des faisceaux coh\'erents sur
\m{C_n}}{dualgen}

Soient \m{P\in C} et $M$ un  \m{\ko_{n,P}}-module de type fini. On note
\m{M^{\vee_n}} le {\em dual} de $M$ :\Nligne
\m{M^{\vee_n}=\Hom(M,\ko_{n,P})}~. Si aucune confusion n'est \`a craindre on
notera \m{M^\vee=M^{\vee_n}}.
Si $N$ est un \m{\ko_{C,P}}-module, on note \m{N^*} le dual de $N$ :
\m{N^*=\Hom(N,\ko_{C,P})}.

Soit $\ke$ un faisceau coh\'erent sur \m{C_n}. On note \m{\ke^{\vee_n}} le {\em
dual} de $\ke$ : \
\m{\ke^{\vee_n}=\HHom(\ke,\ko_n)} . Si aucune confusion n'est \`a craindre on
notera \m{\ke^\vee=\ke^{\vee_n}}.
Si $E$ est un faisceau coh\'erent sur $C$, on note \m{E^*} le dual de $E$ :
\m{E^*=\HHom(E,\ko_C)}.
Ces notations sont justifi\'ees par le fait que \m{E^\vee\not=E^*}. Plus
g\'en\'eralement on a, si $i$ un entier tel que \m{1\leq i\leq n} et $\ke$ un
faisceau coh\'erent sur \m{C_i}, un isomorphisme
canonique
\[\ke^{\vee_n} \ \simeq \ \ke^{\vee_i}\otimes\ki_C^{n-i} ,\]
(\m{\ki_C} d\'esignant le faisceau d'id\'eaux de $C$, qui est un fibr\'e en
droites sur \m{C_{n-1}}). En particulier, pour tout faisceau coh\'erent $E$ sur
$C$, on a \
\m{E^{\vee_n}\simeq E^*\ot L^{n-1}} (cf. \cite{dr4}, lemme 4.1.1) .

Pour tout entier $i$ tel que \m{1\leq i<n}, on a \ \m{(\ke^\vee)^{(i)}=
(\ke_{\mid C_i})^\vee} (\cite{dr4}, prop. 4.1.2).

\sepsubsub 

\begin{subsub}\label{dual} Sous-faisceau de torsion d'un faisceau
coh\'erent sur \m{C_n} - \rm
Soient $P$ un point ferm\'e de $C$ et \m{x\in\ko_{nP}} un \'el\'ement au dessus
d'un g\'en\'erateur de l'id\'eal maximal de \m{\ko_{CP}}.
Soit $M$ un \m{\ko_{nP}}-module de type fini. Le {\em sous-module de torsion
\m{T(M)} de $M$} est constitu\'e des \'el\'ements annul\'es par une puissance de
$x$. On dit que $M$ est {\em sans torsion} si ce sous-module est nul. C'est donc
le cas si et seulement si pour tout \m{m\in M} non nul et tout entier \m{p>0} on
a \m{x^pm\not=0}.

Soit $\ke$ un faisceau coh\'erent sur \m{C_n}. Le {\em sous-faisceau de torsion
\m{T(\ke)} de $\ke$} est le sous-faisceau maximal de $\ke$ dont le support est
fini. Pour tout point ferm\'e $P$ de $C$ on a \ \m{T(\ke)_P=T(\ke_P)}. On a
donc une suite exacte canonique
\[0\lra T(\ke)\lra\ke\lra\ke^{\vee\vee}\lra 0 .\]
\end{subsub}

\sepsubsub

\begin{subsub}\label{refl} Faisceaux r\'eflexifs - \rm Un faisceau coh\'erent
$\ke$ sur \m{C_n} est r\'eflexif si et seulement si il est sans torsion
(\cite{dr4}, th\'eor\`eme 4.2.2), si et seulement si \m{\ke^{(1)}} est
localement libre sur $C$ (\cite{dr4}, prop. 3.3.1).
\end{subsub}

\end{sub}

\sepsub

\Ssect{Dualit\'e des faisceaux de torsion}{du_to}

Soit $\T$ un faisceau de torsion sur \m{C_n}. Alors on a \'evidemment \
\m{\T^\vee=0}. On appelle {\em dual} de $\T$ le faisceau
\[D_n(\T) \ = \ \EExt^1_{\ko_n}(\T,\ko_n) .\]
S'il n'y a pas d'ambigu\" it\'e sur $n$, on notera plus simplement \
\m{\widetilde{\T} = D_n(\T)} .
Rappelons que pour tout \m{i\geq 2} on a \ \m{\EExt^i_{\ko_n}(\T,\ko_n)=0} \
d'apr\`es le corollaire 4.2.4 de \cite{dr4}.

\sepprop

\begin{subsub}\label{tor_pr1}{\bf Proposition : } Soit $\T$ un faisceau de
torsion sur \m{C_i}, \m{1\leq i<n}. Alors on a un isomorphisme canonique
\[D_n(\T) \ \simeq \ D_i(\T)\ot\L^{n-i} .\]
\end{subsub}

Bien s\^ur on a \ \m{D_i(\T)\ot\L^{n-i}\simeq D_i(\T)} .

\begin{proof}
D'apr\`es la proposition 2.3.1 de \cite{dr4} on a un isomorphisme \Nligne
\m{D_i(\T)\simeq\EExt^1_{\ko_n}(\T,\ko_i)}. On consid\`ere la suite exacte
\xmat{0\ar[r] &\ko_i\ot\L^{n-i}\ar[r] & \ko_n\ar[r]^-r & \ko_{n-i}\ar[r] & 0 .}
Il suffit de montrer que le morphisme induit par $r$
\[\Phi:\EExt^1_{\ko_n}(\T,\ko_n)\lra\EExt^1_{\ko_n}(\T,\ko_{n-i})\]
est nul.

On consid\`ere une r\'esolution localement libre de $\T$ :
\xmat{\cdots\E_2\ar[r]^-{f_2} & \E_1\ar[r]^-{f_1} & \E_0\ar[r] & \T\ar[r] & 0 .}
Soient $P$ un point ferm\'e de $C$ et \m{z\in\ko_{nP}} une \'equation de $C$.
Alors \m{\EExt^1_{\ko_n}(\T,\ko_n)} est isomorphe \`a la cohomologie de degr\'e
1 du complexe dual
\xmat{\E_0^\vee\ar[r]^-{^tf_1} & \E_1^\vee\ar[r]^-{^tf_2} & \E_2^\vee\cdots }
et \m{\EExt^1_{\ko_n}(\T,\ko_{n-i})} est isomorphe \`a la cohomologie de degr\'e
1 du complexe obtenu en restreignant le pr\'ec\'edent \`a \m{C_{n-i}}. Le
morphisme $\Phi$ provient du morphisme de complexes
\xmat{\E_0^\vee\ar[r]^-{^tf_1}\flon[d]^{\pi_0} &
\E_1^\vee\ar[r]^-{^tf_2} \flon[d]^{\pi_1} & \E_2^\vee\flon[d]^{\pi_2}\cdots\\
(\E_0^\vee)_{\mid C_{n-i}}\ar[r]^-{^tf_1} &
(\E_1^\vee)_{\mid C_{n-i}}\ar[r]^-{^tf_2} & (\E_2^\vee)_{\mid C_{n-i}}\cdots }
(les fl\`eches verticales \'etant les restrictions). 

Soient $P$ un point du support de $\T$ et \m{z\in\ko_{nP}} une \'equation de
$C$. Soient \ \m{\alpha\in\EExt^1_{\ko_n}(\T,\ko_n)_P} \ et
\m{u\in\ker({}^tf_2)} au dessus de $\alpha$. Puisque $\T$ est concentr\'e sur
\m{C_i}, la multiplication par \m{z^i} : \m{\EExt^1_{\ko_n}(\T,\ko_n)_P\to
\EExt^1_{\ko_n}(\T,\ko_n)_P} \ est nulle. Donc \ \m{z^iu\in\imm({}^tf_1)}, et on
peut \'ecrire \ \m{z^iu={}^tf_1(\theta)}, avec \m{\theta\in(\E^\vee_0)_P}. On va
montrer que $\theta$ est multiple de \m{z^i}. Pour cela on suppose que ce n'est
pas le cas, et on va aboutir \`a une contradiction. On a donc \
\m{\theta=z^k\theta'}, avec \m{0\leq k<i} et \m{\theta'} non multiple de $z$.
On a \ \m{{}^tf_1(z^{n-i+k}\theta')=z^nu=0} , et puisque \m{{}^tf_1} est
injectif, on a \ \m{z^{n-i+k}\theta'=0}. Puisque \m{n-i+k<n}, il en d\'ecoule
que $\theta'$ est multiple de $z$, ce qui est la contradiction recherch\'ee.
On peut donc \'ecrire \ \m{\theta=z^i\theta'} , d'o\`u \ \m{
z^i(u-{}^tf_1(\theta'))=0}, et il en d\'ecoule qu'on peut \'ecrire $u$ sous la
forme :
\m{u  =  {}^tf_1(\theta')+z^{n-i}\rho} .
Il en d\'ecoule que \
\m{\pi_1(u) = {}^tf_1(\pi_0(\theta'))} .
On a donc \ \m{\Phi_P(\alpha)=0}.
\end{proof}

\sepprop

\begin{subsub}\label{tor_co1}{\bf Corollaire : } Soit $\T$ un faisceau de
torsion sur \m{C_n}. Alors on a \ \m{h^0(\T)=h^0(\widetilde{\T})} .
\end{subsub}

\begin{proof}
D'apr\`es la proposition \ref{tor_pr1}, on a, pour tout faisceau de torsion $T$
sur $C$, \m{D_n(T)\simeq T}. Le corollaire en d\'ecoule, en utilisant par
exemple la premi\`ere filtration canonique de $\T$.
\end{proof}

\sepprop

Les faisceaux de torsion sur \m{C_n} et les morphismes entre eux constituent une
cat\'egorie ab\'elienne et no\'eth\'erienne \m{\kt_n(C_n)}, qui est \'evidemment
une sous-cat\'egorie pleine de celle des faisceaux coh\'erents sur \m{C_n}. La
dualit\'e d\'efinit un foncteur contravariant exact
\[D_n:\kt_n(C_n)\lra\kt_n(C_n) .\]

\sepprop

\begin{subsub}\label{tor_pr2}{\bf Proposition : } Le foncteur \m{D_n} est une
involution. Donc si $\T$ est un faisceau de torsion sur \m{C_n}, il
existe un isomorphisme canonique \ \m{\widetilde{\widetilde{\T}}\simeq \T} . 
\end{subsub}

\begin{proof}
Il existe un fibr\'e vectoriel $\E$ et un morphisme surjectif \ \m{f:\E\to\T} .
Alors \m{\ke=\ker(f)} est un faisceau sans torsion, donc r\'eflexif. On obtient
donc en dualisant la suite exacte \ \m{0\to\ke\to\E\to\T\to 0} \ les suivantes
\[0\lra\E^\vee\lra\ke^\vee\lra\widetilde{\T}\lra 0 , \quad\quad
0\lra\ke\lra\E\lra\widetilde{\widetilde{\T}}\to 0\]
d'o\`u on d\'eduit ais\'ement le r\'esultat.
\end{proof}

\sepprop

Si $\T$ est un faisceau de torsion sur \m{C_n}, l'entier \m{h^0(\T)} s'appelle
la {\em longueur} de $T$. On a
\[h^0(\T) \ = \ \sigg_{P\in C}\dim_\C(\T_P) .\]

\sepprop

\begin{subsub}\label{tor_pr3}{\bf Lemme :} Soit $\T$ un faisceau de
torsion sur \m{C_n}. Alors on a  \m{h^0(G_i(\T))=h^0(G^{(i+1)}(\T))} \
pour \m{0\leq i<n}.
\end{subsub}

\begin{proof}
D\'ecoule ais\'ement du corollaire 3.1.8 de \cite{dr4}.
\end{proof}

\sepprop

\begin{subsub}\label{tor_co4}{\bf Corollaire :} Soit $\T$ un faisceau de
torsion sur \m{C_n}. Alors on a, pour \m{1\leq i\leq n} des isomorphismes
canoniques
\[[\widetilde{\T}]_i\simeq\widetilde{[\T_i]}\ot\L^i \ , \quad\quad
(\widetilde{\T})^{(i)}\simeq\widetilde{\T/\T_i}, \quad\quad
G^{(i+1)}(\widetilde{\T})\simeq \widetilde{G_i(\T)} .\]
\end{subsub}

\begin{proof} De la suite exacte \ \m{0\to\T_i\to\T\to\T/\T_i\to 0} \ on
d\'eduit la suivante :
\[0\lra\widetilde{\T/\T_i}\lra\widetilde{\T}\lra\widetilde{[\T_i]}\lra 0 .\]
D'apr\`es la proposition \ref{tor_pr1}, \m{\widetilde{\T/\T_i}} est concentr\'e
sur \m{C_i}. On a donc \ \m{\widetilde{\T/\T_i}\subset(\widetilde{\T})^{(i)}}.
Mais le lemme \ref{tor_pr3} entraine que \ \m{h^0(\widetilde{\T/\T_i})=
h^0((\widetilde{\T})^{(i)})}, donc on a en fait l'\'egalit\'e. Il en d\'ecoule
que \ \m{\widetilde{[\T_i]}\simeq[\widetilde{\T}]_i\ot\L^{-i}} .

Le dernier isomorphisme d\'ecoule de la suite exacte
\[0\lra G^{(i+1)}(\widetilde{\T})\lra[\widetilde{\T}]_i\ot\L\lra
[\widetilde{\T}]_{i+1}\lra 0\]
(cf. lemme 3.1.6 de \cite{dr4}), du fait que par d\'efinition on a
\m{G_i(\T)=\T_i/\T_{i+1}}, et du premier isomorphisme.
\end{proof}

\end{sub}

\sepsub

\Ssect{Dualit\'e des faisceaux sans torsion}{fil_dua}

Soit $\ke$ un faisceau coh\'erent sans torsion sur \m{C_n}. Il est donc
r\'eflexif (cf. \ref{refl}). Les faisceaux \m{\ke_i}, \m{\ke^{(i)}} le sont
donc aussi, \'etant des sous-faisceaux de $\ke$.
Mais les faisceaux \m{\ke/\ke_i} ne le sont pas en g\'en\'eral. On note
\m{\Sigma_i(\ke)} le sous-faisceau de torsion de \m{\ke/\ke_i}, et \m{T_i(\ke)}
celui de \m{G_i(\ke)}.

Pour \m{1\leq i<n}, on note \m{\ke[i]} le noyau du morphisme canonique surjectif
\xmat{\ke\flon[r] & \ke_{\mid C_i}\flon[r] & (\ke_{\mid C_i})^{\vee\vee} .}

\sepsub

\begin{subsub}\label{tor_pr5}{\bf Proposition : } Soit $\ke$ un faisceau
coh\'erent sans torsion sur \m{C_n}. Alors, pour \m{1\leq i<n},

1 - On a un isomorphisme \
\m{\Sigma_i(\ke^\vee) \simeq \widetilde{\Sigma_i(\ke)}\ot\L^i} ,
et une suite exacte
\[0\lra(\ke^\vee)_i\lra(\ke_i)^\vee\ot\L^i\lra\Sigma_i(\ke^\vee)\lra 0 \]
canoniques.

2 - On a un isomorphisme canonique \
\m{\ke[i]^\vee\simeq(\ke^\vee)_i\ot\L^{-i}}.

3 - Il existe un morphisme canonique \
\m{\phi_i(\ke):\Sigma_{i+1}(\ke)\to\Sigma_i(\ke)} \
tel que \ \m{\ker(\phi_i(\ke))\simeq T_i(\ke)} , et que \
\m{\coker(\phi_i(\ke))=R_i(\ke)} \ soit concentr\'e sur $C$.

4 - Il existe une inclusion canonique
\xmat{G^{(i+1)}(\ke^\vee)\flinc[r] & G_i(\ke)^*\ot L^{n-1}}
telle que le quotient soit isomorphe \`a \m{R_i(\ke)}.
\end{subsub}

\begin{proof}
En dualisant la suite exacte \ \m{0\to\ke_i\to\ke\to\ke/\ke_i\to 0} , on
obtient la suite exacte
\[0\lra(\ke/\ke_i)^\vee\lra\ke^\vee\lra(\ke_i)^\vee\lra
 \EExt^1_{\ko_n}(\ke/\ke_i,\ko_n)=\widetilde{\Sigma_i(\ke)}\lra 0 .\]
D'apr\`es la proposition 4.1.2 de \cite{dr4} on a \
\m{(\ke/\ke_i)^\vee=(\ke^\vee)^{(i)}}. On en d\'eduit la suite exacte
\begin{equation}\label{eq3}
0\lra(\ke^\vee)_i\ot\L^{-i}\lra(\ke_i)^\vee\lra\widetilde{\Sigma_i(\ke)}\lra 0 .
\end{equation}
En la dualisant et tensorisant par \m{\L^{-i}}, et en utilisant la proposition
 \ref{tor_pr2} on obtient la suite exacte suivante :
\begin{equation}\label{eq4}
0\lra\ke_i\ot\L^{-i}\lra((\ke^\vee)_i)^\vee\lra\Sigma_i(\ke)\ot\L^{-i}\lra 0 ,
\end{equation}
qui est (\ref{eq3}) avec \m{\ke^\vee} \`a la place de $\ke$. On obtient donc
l'isomorphisme canonique de 1-. On en d\'eduit 2- en dualisant la suite exacte
\ \m{0\to\ke_i\to\ke[i]\to\Sigma_i(\ke)\to 0} .

On a un diagramme commutatif avec lignes et colonnes exactes
\xmat{
        & 0\ar[d] & 0\ar[d]\\
0\ar[r] & \ke_{i+1}\ar[r]\ar[d] & ((\ke^\vee)_{i+1})^\vee\ot\L^{i+1}\ar[r]\ar[d]
& \Sigma_{i+1}(\ke)\ar[r] & 0\\
0\ar[r] & \ke_i\ar[r]\ar[d] & ((\ke^\vee)_i)^\vee\ot\L^i\ar[r]\ar[d] &
\Sigma_i(\ke)\ar[r] & 0\\
 & G_i(\ke)\ar[d] & G^{(i+1)}(\ke^\vee)^\vee\ar[d]\\
 & 0 & 0
}
o\`u les suites horizontales proviennent de (\ref{eq4}) et la suite verticale
du milieu du lemme 3.1.6 de \cite{dr4}. On en d\'eduit ais\'ement 3- et 4-.
\end{proof}
\end{sub}

\sepsub

\Ssect{Invariants du dual d'un faisceau coh\'erent}{invdual}

\begin{subsub}{\bf Proposition : }\label{dua1} Soit $\ke$ un faisceau
coh\'erent. Alors on a
\[R(\ke^\vee)=R(\ke) , \quad\quad
\Deg(\ke^\vee)=-\Deg(\ke)+R(\ke)(n-1)\deg(L)+h^0(T(\ke)) .\]
\end{subsub}

\begin{proof} La premi\`ere assertion concernant les rangs est imm\'ediate, par
exemple en se pla\c cant sur l'ouvert o\`u $\ke$ est quasi localement libre.
D\'emontrons la seconde. Soit \ \m{\kf=\ke/T(\ke)}, qui est un faisceau sans
torsion. On a \ \m{\Deg(\ke)=\Deg(\kf)+h^0(T(\ke))}, \m{R(\ke)=R(\kf)} et
\m{\ke^\vee=\kf^\vee}, donc la seconde assertion \'equivaut \`a
\[\Deg(\kf^\vee)=-\Deg(\kf)+R(\kf)(n-1)\deg(L) .\]
On peut donc supposer que $\ke$ est sans torsion. On va montrer que
\begin{equation}\label{equ}
\Deg(\ke^\vee)=-\Deg(\ke)+R(\ke)(n-1)\deg(L)
\end{equation}
par r\'ecurrence sur $n$. Si \m{n=1} c'est \'evident. Supposons que \m{n>1} et
que (\ref{equ}) soit vraie pour \m{n-1}. On a donc
\[\Deg((\ke_1)^{\vee_{n-1}})=-\Deg(\ke_1)+R(\ke_1)(n-2)\deg(L) .\]
Mais d'apr\`es \ref{dualgen} on a \ \m{(\ke_1)^\vee=(\ke_1)^{\vee_{n-1}}\ot
\ki_C}, donc
\[\Deg((\ke_1)^\vee)=\Deg((\ke_1)^{\vee_{n-1}})+R(\ke_1)\deg(L) ,\]
d'o\`u
\begin{equation}\label{equb}
\Deg((\ke_1)^\vee)=-\Deg(\ke_1)+R(\ke_1)(n-1)\deg(L)
\end{equation}
(c'est-\`a-dire que (\ref{equ}) est vraie pour \m{\ke_1}). D'apr\`es la suite
exacte
\[0\lra\ke_1\lra\ke\lra\ke_{\mid C}\lra 0\]
on a \ \m{\Deg(\ke)=\Deg(\ke_1)+\Deg(\ke_{\mid C})}. Soit \m{T=T(\ke_{\mid C})}.
On a une suite exacte
\[0\lra(\ke_{\mid C})^\vee\lra\ke^\vee\lra(\ke_1)^\vee\lra\widetilde{T}\lra 0,\]
donc
\[\Deg(\ke^\vee)=\Deg((\ke_{\mid C})^\vee)+\Deg((\ke_1)^\vee)-h^0(T) .\]
Mais
\[\Deg((\ke_{\mid C})^\vee)-h^0(T)=-\Deg(\ke_{\mid C})+
(n-1)R(\ke_{\mid C})\deg(L)\]
(car \ \m{(\ke_{\mid C})^\vee=(\ke_{\mid C})^*\ot L^{n-1}}). Donc
\begin{eqnarray*}
\deg(\ke^\vee) & = & \Deg((\ke_1)^\vee)-\Deg(\ke_{\mid C})+
(n-1)R(\ke_{\mid C})\deg(L)\\
 & = & -\Deg(\ke)+R(\ke)(n-1)\deg(L)
\end{eqnarray*}
d'apr\`es (\ref{equb}).
\end{proof}

\sepprop

\begin{subsub}\label{cor1}{\bf Corollaire : } Soit $\ke$ un faisceau
coh\'erent r\'eflexif sur \m{C_n}. Alors, pour \m{1\leq i<n}, on a
\[R\big((\ke^\vee)_i\big)=R(\ke_i) , \quad
R\big((\ke^\vee)^{(i)}\big)=R(\ke^{(i)}), \quad
R\big((\ke^\vee)_{\mid C_i}\big)=R(\ke_{\mid C_i}) ,\]
\[\Deg\big((\ke^\vee)_i\big) \ = \ -\Deg(\ke_i)+(n+i-1)R(\ke_i)\deg(L)
-h^0(\Sigma_i(\ke)) ,\]
\[\Deg\big((\ke^\vee)_{\mid C_i}\big) \ = \ \Deg((\ke_{\mid
C_i})^\vee)-iR(\ke_i)\deg(L) .\]
\end{subsub}

\begin{proof}
D\'ecoule ais\'ement des propositions \ref{tor_pr5} et \ref{dua1}.
\end{proof}

\sepprop

\begin{subsub}\label{cor2}{\bf Corollaire : } Soient $\ke$ un faisceau
coh\'erent r\'eflexif sur \m{C_n} et $i$ un entier tel que \m{1\leq i<n} et
\m{R(\ke_i)>0}. Alors on a
\[\mu\big((\ke^\vee)_{\mid C_i}\big)-\mu\big((\ke^\vee)_i\big) \ = \
\mu(\ke_i\ot\L^{-i})-\mu(\ke^{(i)})+h^0(\Sigma_i(\ke))\big(
\frac{1}{R(\ke^{(i)})}+\frac{1}{R(\ke_i)}\big) .\]
\end{subsub}

\end{sub}

%\sepsec
\newpage

\section{Conditions d'existence des faisceaux (semi-)stables}\label{condex}

Dans toute la suite de ce chapitre on consid\`ere une courbe multiple primitive
\m{C_n} de courbe r\'eduite associ\'ee $C$. On utilise les notations de
\ref{def_nota}, et on suppose que \m{\deg(L)<0}.

\sepsub

\Ssect{Crit\`eres de (semi-)stabilit\'e}{critstab}

\begin{subsub}\label{lem1}{\bf Lemme : } Soient $A$, $A''$, $B$, $B''$, $E$,
$E''$ des faisceaux coh\'erents de rang positif sur \m{C_n}, tels que
\[R(E)=R(A)+R(B) , \quad R(E'')=R(A'')+R(B'') ,\]
\[\Deg(E)=\Deg(A)+\Deg(B) , \quad \Deg(E'')=\Deg(A'')+\Deg(B'') .\]
On suppose qu'on a \ \m{\mu(B)\geq\mu(A)}, \m{\mu(A'')\geq\mu(A)},
\m{\mu(B'')\geq\mu(B)}, et que
\[\frac{R(E'')}{R(E)} \ \geq \ \frac{R(A'')}{R(A)} \quad .\]
Alors on a \m{\mu(E'')\geq\mu(E)} . Si de plus \ \m{\mu(A'')>\mu(A)} ou
\m{\mu(B'')>\mu(B)}, alors on a \m{\mu(E'')>\mu(E)} .
\end{subsub}

\begin{proof}
D'apr\'es les hypoth\'eses
\[\frac{R(E'')}{R(E)}\geq\frac{R(A'')}{R(A)}\]
\'equivaut \`a
\[\frac{R(B'')}{R(B)}\geq\frac{R(A'')}{R(A)} ,\]
et
\[\mu(E'')-\mu(E) \ = \ \frac{\Delta}{R(E)R(E'')} ,\]
avec
\[
\Delta=\big(\Deg(A'')+\Deg(B'')\big)\big(R(A)+R(B)\big)-
\big(\Deg(A)+\Deg(B)\big)\big(R(A'')+R(B'')\big) .\]
On a
\[\Deg(A'')\geq\Deg(A)\frac{R(A'')}{R(A)} , \quad\quad
\Deg(B'')\geq\Deg(B)\frac{R(B'')}{R(B)} ,\]
Donc \m{\Delta\geq\Delta'}, avec
\begin{eqnarray*}
\Delta' & = & \big(\Deg(A)\frac{R(A'')}{R(A)}+\Deg(B)\frac{R(B'')}{R(B)}\big)
\big(R(A)+R(B)\big)-\\
& & \quad\quad\quad\quad\big(\Deg(A)+\Deg(B)\big)\big(R(A'')+R(B'')\big)\\
 & = & \big(\mu(B)-\mu(A)\big)\big(R(B'')R(A)-R(A'')R(B)\big) .
\end{eqnarray*}
Le r\'esultat en d\'ecoule imm\'ediatement.
\end{proof}

\sepprop

\begin{subsub}\label{theo1}{\bf Th\'eor\`eme : } Soient $\ke$ un faisceau
coh\'erent sans torsion sur \m{C_n} et $k$ un entier tel que \m{1\leq k<n} et
que \m{\ke_k\not=0}. On suppose que
\begin{equation}\label{eqX}\mu(\ke^{(k)})\leq\mu(\ke) ,\quad
\mu((\ke^\vee)^{(k)})\leq\mu(\ke^\vee) .\end{equation}
Alors, si \m{\ke[k]}, \m{(\ke_{\mid C_k})^{\vee\vee}}, \m{(\ke^\vee)[k]} et
\m{((\ke^\vee)_{\mid C_k})^{\vee\vee}} sont semi-stables il en est de m\^eme de
$\ke$.

Si de plus les in\'egalit\'es de (\ref{eqX}) sont strictes, et si \m{\ke[k]} ou
\m{(\ke_{\mid C_k})^{\vee\vee}}, ainsi que \m{(\ke^\vee)[k]} ou
\m{((\ke^\vee)_{\mid C_k})^{\vee\vee}} sont stables, alors $\ke$ est stable.
\end{subsub}

\begin{proof}
Supposons que les hypoth\`eses du th\'eor\`eme soient v\'erifi\'ees. Soit 
\xmat{\ke\flon[r] & \ke''}
 un quotient de $\ke$. Il faut montrer que \m{\mu(\ke'')\geq\mu(\ke)} .
On peut supposer que  \m{\ke''} est sans torsion. On a un diagramme commutatif
avec lignes exactes et fl\`eches verticales surjectives
\xmat{0\ar[r] & \ke[k]\ar[r]\flon[d] & \ke\ar[r]\flon[d] & 
(\ke_{\mid C_k})^{\vee\vee}\ar[r]\flon[d] & 0\\
0\ar[r] & \ke''[k]\ar[r] & \ke''\ar[r] & (\ke''_{\mid C_k})^{\vee\vee}\ar[r] & 0
}
Le cas o\`u \m{\ke''[k]=0} est \'evident. On supposera donc que
\m{\ke''[k]\not=0} . Remarquons que les in\'egalit\'es (\ref{eqX}) \'equivalent
\`a
\[\mu((\ke_{\mid C_k})^{\vee\vee})\geq\mu(\ke[k]) ,\quad
\mu(((\ke^\vee)_{\mid C_k})^{\vee\vee})\geq\mu(\ke^\vee[k])\]
(car \ \m{(\ke^\vee)^{(k)}=(\ke_{\mid C_k})^\vee}).
Le morphisme vertical de droite du diagramme pr\'ec\'edent est surjectif, donc
on a \ \m{\mu((\ke''_{\mid C_k})^{\vee\vee})\geq
\mu((\ke_{\mid C_k})^{\vee\vee})} \ d'apr\`es la semi-stabilit\'e de
\m{(\ke_{\mid C_k})^{\vee\vee}}. Le conoyau du morphisme vertical de gauche est
de torsion, donc on a \ \m{\mu(\ke''[k])\geq\mu(\ke[k])} \ d'apr\`es la
semi-stabilit\'e de \m{\ke[k]}.
D'apr\`es le lemme \ref{lem1} on a \m{\mu(\ke'')\geq\mu(\ke)} si
\[\frac{R(\ke'')}{R(\ke)} \ \geq \ \frac{R(\ke''[k])}{R(\ke[k])}\quad .\]
On peut donc supposer que
\begin{equation}\label{equ3}
\frac{R(\ke'')}{R(\ke)} \ < \ \frac{R(\ke''[k])}{R(\ke[k])}\quad .
\end{equation}
On utilise maintenant la suite exacte
\[0\lra\ke''^\vee\lra\ke^\vee\lra\ke'^\vee\lra 0 \]
obtenue en utilisant le fait que \m{\ke''} est r\'eflexif.
D'apr\`es le lemme \ref{dua1}, \m{\mu(\ke'')\geq\mu(\ke)} \'equivaut \`a
\m{\mu(\ke'^\vee)\geq\mu(\ke^\vee)}, et d'apr\`es le lemme \ref{lem1}, cette
in\'egalit\'e est v\'erifi\'ee si
\begin{equation}\label{equ1}
\frac{R(\ke')}{R(\ke)}=\frac{R(\ke'^\vee)}{R(\ke^\vee)} \ \geq \
\frac{R(\ke'^\vee[k])}{R(\ke^\vee[k])} \quad .
\end{equation}
D'apr\`es \ref{dual}, on a \ \m{R(\ke'^\vee[k])=R(\ke'[k])} \ et \
\m{R(\ke^\vee[k])=R(\ke[k])} . Donc (\ref{equ1}) \'equivaut \`a
\begin{equation}\label{equ2}
\frac{R(\ke')}{R(\ke)} \ \geq \ \frac{R(\ke'[k])}{R(\ke[k])} \quad .
\end{equation}
Puisque \m{\ke'_k} est contenu dans le noyau du morphisme canonique surjectif \
\m{\ke_k\to\ke''_k} , on a
\[R(\ke'[k])=R(\ke'_k) \ \leq \ R(\ke_k)-R(\ke''_k)=R(\ke[k])-R(\ke''[k]) ,\]
donc on
peut \'ecrire
\[R(\ke'[k])=R(\ke[k])-R(\ke''[k])-\eta ,\]
avec \m{\eta\geq 0}. Donc (\ref{equ2}) s\'ecrit
\[\frac{R(\ke'')}{R(\ke)} \ \leq \ \frac{R(\ke''[k])+\eta}{R(\ke[k])}\quad .\]
L'in\'egalit\'e pr\'ec\'edente est vraie d'apr\`es (\ref{equ3}). On a donc bien
\ \m{\mu(\ke'')\geq\mu(\ke)} .

L'assertion concernant la stabilit\'e se d\'emontre de mani\`ere analogue.
\end{proof}

\end{sub}

\sepsub

\Ssect{Le cas des fibr\'es vectoriels}{vect}

\begin{subsub}\label{theo2}{\bf Th\'eor\`eme : } Soit $\E$ un fibr\'e vectoriel
sur \m{C_n}. Alors, si \m{\E_{\mid C}} est semi-stable (resp. stable), il en est
de m\^eme de $\E$.
\end{subsub}

\begin{proof}
Posons \m{E=\E_{\mid C}}. Alors les filtrations canoniques de $\E$ sont
identiques, et leurs gradu\'es sont
\[\big(G_0(\ke),G_1(\ke),\ldots,G_{n-1}(\ke)\big) \ = \ (E,E\ot L,\ldots,E\ot
L^{n-1}) \ .\]
Les in\'egalit\'es (\ref{eqX}) sont trivialement v\'erifi\'ees (car
\m{\deg(L)<0}) pour tout entier $k$ tel que \m{1\leq k<n}.

Le th\'eor\`eme \ref{theo2} se d\'emontre par r\'ecurrence sur $n$ : pour
\m{n=1} c'est \'evident. Supposons que ce soit vrai pour \m{n-1\geq 1}. Alors
\m{\E_1} est semi-stable (resp. stable). Le th\'eor\`eme \ref{theo1} permet
alors de conclure que $\E$ est semi-stable (resp. stable).
\end{proof}

\sepprop

\begin{subsub} Vari\'et\'es de modules - \rm Soient $r$, $\delta$ des entiers
tels que \m{r\geq 1}. Posons
\[R = nr , \quad d = n\delta+\frac{n(n-1)}{2}r\deg(L) .\]
Pour tout fibr\'e vectoriel $\E$ sur \m{C_n} tel que \m{\E_{\mid C}} soit de
rang $r$ et de degr\'e $\delta$, on a \ \m{R(\E)=R} \ et \ \m{\Deg(\ke)=d} . Il
d\'ecoule de \ref{constcoh3} que la vari\'et\'e de modules \m{\M(R,d)} des
fibr\'es vectoriels stables de rang g\'en\'eralis\'e $R$ et de degr\'e
g\'en\'eralis\'e $d$ sur \m{C_n} est non vide. C'est un ouvert irr\'eductible et
lisse de la vari\'et\'e de modules \m{\km(R,d)} des faisceaux stables de rang
g\'en\'eralis\'e $R$ et de degr\'e g\'en\'eralis\'e $d$.
\end{subsub}

\end{sub}

\sepsub

\Ssect{Le cas des faisceaux quasi localement libres de type g\'en\'erique}{qllx}

\begin{subsub}\label{theo3}{\bf Th\'eor\`eme : } Soient $a$, $k$ des entiers
tels que \m{a>0} et \m{1\leq k<n}. Soit $\ke$ un faisceau quasi localement libre
de type rigide, localement isomorphe \`a \m{a\ko_n\oplus\ko_k} et tel que
\begin{equation}\label{equCC3}
\mu(V_\ke)+\frac{n}{2}\deg(L)\leq\mu(F_\ke)\leq\mu(E_\ke)-\frac{n}{2}\deg(L).
\end{equation}

Alors si \m{E_\ke}, \m{F_\ke} et \m{V_\ke} sont semi-stables, il en est de
m\^eme de $\ke$.

Si les in\'egalit\'es pr\'ec\'edentes sont strictes, et si \m{E_\ke}, \m{F_\ke}
et \m{V_\ke} sont stables, il en est de m\`eme de $\ke$.
\end{subsub}

\begin{proof} On ne d\'emontrera que la premi\`ere assertion, la seconde
\'etant analogue. On utilise les notations de \ref{def_qlltr}. Supposons les
in\'egalit\'es (\ref{equCC3}) v\'erifi\'ees et \m{E_\ke}, \m{F_\ke} et
\m{V_\ke} semi-stables. Alors on a
\[\ke[k]=\ke_k=\F\ ,\quad\quad \ke_{\mid C_k}=\E ,\quad\quad \ke^\vee[k]=
(\ke^\vee)_k=\F\ot\L^k\ ,
\quad\quad (\ke^\vee)_{\mid C_k}=\V^\vee\  .\]
Donc d'apr\`es le th\'eor\`eme \ref{theo2}, \m{\ke[k]}, \m{\ke_{\mid C_k}},
\m{\ke^\vee[k]} et \m{(\ke^\vee)_{\mid C_k}} sont semi-stables. Un calcul
simple montre que les in\'egalit\'es (\ref{equCC3}) \'equivalent aux
in\'egalit\'es (\ref{eqX}). La semi-stabilit\'e de $\ke$ d\'ecoule donc du
th\'eor\`eme \ref{theo1}.
\end{proof}

\sepprop

La semi-stabilit\'e de \m{E_\ke}, \m{F_\ke} et \m{V_\ke} entrainent d'autres
in\'egalit\'es :
\[\mu(E_\ke)\leq\mu(F_\ke) , \quad\quad \mu(F_\ke)\leq\mu(V_\ke) \]
(car il existe un morphisme surjectif \m{E_\ke\to F_\ke} et un morphisme
injectif \m{F_\ke\to V_\ke}). Les in\'egalit\'es pr\'ec\'edentes et
(\ref{equCC3}) \'equivalent aux in\'egalit\'es
\[\mu(E_\ke) \ \leq\mu(F_\ke) \ \leq\mu(E_\ke)-\frac{n-k}{a+1}\deg(L) .\]

\sepprop

\begin{subsub}Vari\'et\'es de modules de faisceaux stables - \rm Soient $a$,
$k$, $\epsilon$, $\delta$ des entiers, avec \m{a\geq 1} et \m{1\leq k<n}. Soient
\[R \ = \ an+k , \quad\quad\quad d \ = \ k\epsilon+(n-k)\delta+(n(n-1)a+
k(k-1))\frac{\deg(L)}{2} \ .\]
On note \m{\km(R,d)} la vari\'et\'e de modules des faisceaux stables de rang
g\'en\'eralis\'e $R$ et de degr\'e g\'en\'eralis\'e $d$ sur \m{C_n}. Les
faisceaux quasi localement libres $\ke$ de type g\'en\'erique stables localement
isomorphes \`a \m{a\ko_n\oplus\ko_k} et tels que \m{E_\ke} (resp. \m{F_\ke})
soit de rang \m{a+1} (resp. $a$) et de degr\'e $\epsilon$ (resp. $\delta$)
constituent un ouvert irr\'eductible de \m{\km(R,d)}, dont la
sous-vari\'et\'e r\'eduite associ\'ee est not\'ee \m{\kn(a,k,\delta,\epsilon)}.
\end{subsub}

\sepprop

\begin{subsub}\label{theo4}{\bf Th\'eor\`eme : } Si on a
\[\frac{\epsilon}{a+1} \ < \ \frac{\delta}{a} \ < \
\frac{\epsilon-(n-k)\deg(L)}{a+1}\]
\m{\kn(a,k,\delta,\epsilon)} est non vide.
\end{subsub}

\begin{proof}Les hypoth\`eses et \cite{rute} impliquent qu'il existe des
fibr\'es stables $E$, $F$, $V$ sur $C$, tels que
\[rg(E)=a+1, \quad \deg(E)=\epsilon, \quad \rg(F)=a, \quad \deg(F)=\delta,\]
\[\quad \rg(V)=a+1, \quad \deg(V)=\epsilon-(n-k)\deg(L)\]
tels qu'il existe une suite exacte
\[0\lra F\ot L^{n-k}\lra V\ot L^{n-k}\lra E\lra F\lra 0 .\]
D'apr\`es la proposition \ref{propC3} il existe un faisceau quasi localement
libre de type rigide $\ke$ , localement isomorphe \`a \m{a\ko_n\oplus\ko_k} et
tel que \m{(*)_\ke} soit isomorphe \`a la suite exacte pr\'ec\'edente.
D'apr\`es le th\'eor\`eme \ref{theo3}, $\ke$ est stable, et d\'efinit donc un
point de \m{\kn(a,k,\delta,\epsilon)}.
\end{proof}

\end{sub}

\sepsub

\Ssect{Exemple d'application \`a des faisceaux non quasi localement libres}{ex2}

Soient $\E$ un fibr\'e vectoriel sur \m{C_n}, \m{E=\E_{\mid C}} et $Z$ un
ensemble fini de points de $C$. On pose \m{z=h^0(\ko_Z)}. Soient
\m{\phi:\E\to\ko_Z} un morphisme surjectif, et \ \m{\ke_\phi=\ker(\phi)}. On a
deux suites exactes
\[0\lra\ke_\phi\lra\E\lra\ko_Z\lra 0 , \quad
0\lra\E^\vee\lra\ke_\phi^\vee\lra\ko_Z\lra 0 .\]
Le morphisme $\phi$ se factorise par $E$. On note \m{E_\phi} le noyau du
morphisme induit \m{E\to\ko_Z}. On note \m{\ke'_\phi} le noyau du morphisme
induit \ \m{\E_{\mid C_{n-1}}\to\ko_Z}.

\sepprop

\begin{subsub}\label{lemX}{\bf Lemme : } On a \ \m{\ke_\phi[1]=\E_1} ,
\m{(\ke_{\phi\mid C})^{\vee\vee}=E_\phi} ,
\m{\ke_\phi^\vee[1]=(\ke'_\phi)^\vee} \ et \Nligne
\m{\big((\ke_\phi)^\vee)_{\mid C}\big)^{\vee\vee}=E^*}.
\end{subsub}

\begin{proof} Il suffit de le d\'emontrer en un point $P$ de $Z$. Soit
\m{z\in\ko_{n,P}} une \'equation de $C$ et \m{x\in\ko_{n,P}} au dessus d'un
g\'en\'erateur de l'id\'eal maximal de $P$ dans \m{\ko_C}. Si \m{r=\rg(\E_{\mid
C})}, on a \ \m{\ke_{\phi,P}\simeq r\ko_{n,P}\oplus(x,z)}. On peut donc supposer
que \m{\ke_{\phi,P}=(x,z)}. Il faut montrer que \ \m{\ke_\phi[1]_P=(z)}.
On a \ \m{(x,z)_{\mid C}=(x)/(xz)\oplus(z)/(xz,z^2)}. Le
premier facteur est isomorphe \`a \m{\ko_{C,P}} et le second \`a $\C$. Donc
\m{\ke_\phi[1]_P} est le noyau du morphisme
\[{\xymatrix@R=6pt{
(x,z)\ar[r] & \ko_{C,P} \\
\ \ \ \ ax+bz\fmaps[r] & \ov{a} }}\]
($\ov{a}$ d\'esignant l'image de $a$ dans \m{\ko_{C,P}}). On a donc \
\m{\ke_\phi[1]_P=(z)=\E_{1,P}}. On a \ \m{\ke_{\phi\mid C}=E_\phi\oplus\ko_Z},
donc \ \m{(\ke_{\phi\mid C})^{\vee\vee}=E_\phi}.

On a \ \m{\ke_\phi^\vee[1]=((\ke_\phi)_1)^\vee\ot\L} d'apr\`es la proposition
\ref{tor_pr5}, 2-. On a un diagramme commutatif avec lignes exactes
\xmat{0\ar[r] & E\ot L^{n-1}=(\ke_\phi)^{(1)}\ar[r]\fleq[d] &
\ke_\phi\ar[r]\flinc[d] & (\ke_\phi)_1\ot\L^{-1}\ar[r]\flinc[d] & 0\\
0\ar[r] & E\ot L^{n-1}=\E^{(1)}\ar[r] & \E\ar[r] & \E_{\mid C_{n-1}}\ar[r] & 0}
On en d\'eduit imm\'ediatement la troisi\`eme \'egalit\'e. On a enfin
\[((\ke^\vee)_{\mid C})^{\vee\vee}=(\ke^{(1)})^\vee=(E\ot L^{n-1})^\vee=E^* .\]
\end{proof}

\sepprop

\begin{subsub}\label{theo5}{\bf Th\'eor\`eme : } Si on a \ \m{z\leq
-rg(E)\deg(L)} (resp. $<$) et si $E$ et \m{E_\phi} sont semi-stables (resp.
stables), alors $\ke_\phi$ est semi-stable (resp. stable).
\end{subsub}

\begin{proof} Cela se d\'emontre ais\'ement par r\'ecurrence sur $n$, en
utilisant le lemme \ref{lemX} et les th\'eor\`emes \ref{theo1} et \ref{theo2}.
\end{proof}

\end{sub}

\end{document}